\documentclass[draftclsnofoot,11pt,onecolumn]{IEEEtran}
\usepackage{amsmath}

\usepackage{amsmath,bm}
\usepackage{amsfonts}
\usepackage{amssymb}
\usepackage{graphicx}
\usepackage{tabularx}
\usepackage[left=2cm,right=2cm,top=2cm,bottom=2cm]{geometry}
\usepackage{hyperref}
\usepackage{algorithm} 
\usepackage{savesym}
\savesymbol{AND}
\savesymbol{OR}
\savesymbol{NOT}
\savesymbol{TO}
\savesymbol{COMMENT}
\savesymbol{BODY}
\savesymbol{IF}
\savesymbol{ELSE}
\savesymbol{ELSIF}
\savesymbol{FOR}
\savesymbol{WHILE}
\usepackage{algorithmic}

\usepackage{amsthm}
\newtheorem{Definition}{Definition}
\newtheorem{Theorem}{Theorem}
\newtheorem{Lemma}{Lemma}
\newtheorem{Remark}{Remark}

\newtheorem{Corollary}{Corollary}

\newtheorem{Example}{Example}
\usepackage{color}
\usepackage{caption}
\usepackage{subcaption}
\usepackage{cite}
\usepackage{url}

\usepackage{color}

\newcommand{\m}[1]{\mathbf{#1}}
\newcommand{\mc}[1]{\mathcal{#1}}
\newcommand{\mb}[1]{\mathbb{#1}}

\begin{document}

\title{Theory and Applications of Matrix-Weighted Consensus}
\author{Minh~Hoang~Trinh$^{\dag}$ 
and Hyo-Sung~Ahn$^{\dag}$
\thanks{\small $^{\dag}$School of Mechanical Eng., Gwangju Institute of Science and Technology, Gwangju, Korea.
{E-mails: \{trinhhoangminh,hyosung\}@gist.ac.kr}}
}

%

\markboth{Technical Report}{}%
%


\maketitle

\begin{abstract}
This paper proposes the \emph{matrix-weighted consensus algorithm}, which is a generalization of the consensus algorithm in the literature. Given a networked dynamical system where the interconnections between agents are weighted by nonnegative definite matrices instead of nonnegative scalars, consensus and clustering phenomena naturally exist. We examine algebraic and algebraic graph conditions for achieving  a consensus, and provide an algorithm for finding all clusters of a given system. Finally, we illustrate two applications of the proposed consensus algorithm in cluster consensus and in bearing-based formation control.
\end{abstract}


\ifCLASSOPTIONpeerreview
\begin{center} \bfseries EDICS Category: 3-BBND \end{center}
\fi

\IEEEpeerreviewmaketitle

\section{Introduction}
\label{sec:1}
\IEEEPARstart{C}{onsensus} algorithm has been extensively studied in the literature as a main tool for solving the cooperative control problems in multiagent systems \cite{Olfati2004,Olfati2007,ren2005survey}. In fact, consensus algorithm and its modifications are found in broad applications, for examples, in control of unmanned vehicle formations \cite{Fax2004,Lafferriere2005,Ren2007scl,Oh2015survey}, network synchronization \cite{Kim2013aut,Tuna2009tac}, modeling social networks \cite{Hendrickx2014,Xia2016}, and coordination of power distribution systems and automated traffic networks \cite{Kim2016aut,Kim2016sj}.

Given a system of $n$ single-integrator agents whose interconnections between agents are characterized by an weighted undirected graph $\mc{G}$, the consensus algorithm \cite{Olfati2004} is defined as\footnote{Formal definitions will be provided in the next section.}
\begin{equation} \label{eq:consensus-olfati}
\dot{\m{x}}_i = \sum_{j =1}^n a_{ij}(\m{x}_j - \m{x}_i), \forall i =1, \ldots, n,
\end{equation}
where $\m{x}_i,~\m{x}_j \in \mb{R}^d$ are the state vectors of agents $i$ and $j$, and $a_{ij}$ is a positive scalar (or zero) if $i$ and $j$ are connected (or disconnected, respectively). It is well-known that under the consensus protocol \eqref{eq:consensus-olfati}, an average consensus is globally achieved if and only if $\mc{G}$ is connected \cite{Olfati2004}.

This paper generalizes the consensus algorithm \eqref{eq:consensus-olfati} by using matrix weight $\m{A}_{ij}$ instead of the scalar weight $a_{ij}$ to describe the interconnection between two agents $i$ and $j$. Here, a matrix weight could be a positive definite matrix (strong connection), a positive semidefinite matrix (weak connection), or a zero matrix (no direct connection). Thus, the matrix-weight consensus covers a larger set of problems in multi-agent systems. 

In the literature, matrix weights arise in many problems to describe the interconnections between agents. For example, the author of \cite{Tuna2016aut} used matrix weights to describe interconnections between coupled linear oscillators and provided conditions to synchronize these networks in some situations. The concept of deviated cyclic pursuit introduced in \cite{ramirez2009distributed}, and orientation estimation in \cite{lee2016aut,ahn2017consensus} can be considered as consensus protocols with rotation matrix weights. Also, the bearing-based formation control setup in \cite{Zhao2015CDC} can be formulated as a special case of the matrix-weighted consensus protocol proposed in this paper. In the context of social networks, suppose that a group of people are discussing multiple topics, matrix weights were used to describe the logical inter-dependency of the topics \cite{friedkin2016network,parsegov2017novel}. However, the works \cite{friedkin2016network,parsegov2017novel} only considered a discrete-time model in which the matrix weights are the same for all edges.

In this paper, we study the matrix-weight consensus algorithm with  undirected graphs. We firstly define several terminologies (for e.g., positive/semipositive connections, positive tree, matrix-weighted Laplacian, etc), and prove some basic algebraic properties of the matrix-weighted graph. Secondly, we propose the matrix-weight consensus protocol and provide a necessary and sufficient condition for globally exponentially reaching an average consensus based on the nullspace of the matrix-weighted Laplacian. Next, due to the existence of semidefinite matrix weights, clustered consensus happens naturally even when the graph is connected. We examine the algebraic graph theory of consensus and clustering phenomena. Further, an algorithm to determine all clusters in the network is provided. The algorithm initially partitions the graph into a set of clusters associated with the positive trees in the graph. If two clusters satisfy several algebraic conditions on their connections, they will be merged together at each iteration of the algorithm. The algorithm gradually reduces the number of clusters in the graph, and it ends when no two  clusters can be further merged together. If there is a cluster containing all vertices in the graph, under matrix-weighted consensus protocol, a consensus is globally achieved. Otherwise, we know the exact number of clusters in the system under the matrix-weight consensus protocol. Finally, two examples are given to illustrate applications of the proposed matrix-weight consensus algorithm. The first example demonstrates how clustered consensus can be used to gather a group of agents into several clusters. The second example is taken from the bearing-based formation control in the literature \cite{Zhao2015CDC,Zhao2015CNS}.

The rest of this paper is organized as follows. Section \ref{sec:2} defines basic terminologies and introduces the matrix-weighted consensus algorithm. The algebraic condition for globally reaching a consensus in undirected networks is presented in Section \ref{sec:3}. Section \ref{sec:4} further studies the consensus and clustered consensus phenomena under the matrix-weighted consensus algorithm. Section \ref{sec:5} provides two applications of the proposed algorithm. Finally, Section \ref{sec:6} summarizes the paper and discusses several further research directions.

\subsection{Notations}
In this paper, $\mb{R}^d$ denotes the Euclidean $d$-dimensional space. Vectors and matrices are denoted by bold-font letters, while sets are denoted by calligraphic characters. Note $\m{1}_n \in \mb{R}^{n}$ denotes the vector of all entries 1s, $\m{I}_{d}$ denotes the identity matrix in of dimension $d \times d$, and $\otimes$ denotes the Kronecker product.
\section{Preliminaries And Problem Formulation}
\label{sec:2}
\subsection{Matrix-Weighted Graphs}
This subsection sets a framework for introducing the matrix-weighted consensus protocol and the main analysis of this paper. Most of definitions are analogous to the definitions of algebraic graph theory \cite{Godsil2001}.

A fixed undirected graph with matrix weights is denoted by $\mc{G}$. The graph $\mc{G}$ is characterized by a triple $(\mc{V}, \mc{E}, \mc{A})$. Here, $\mc{V} = \{ 1, \ldots, n \}$ denotes the set of $|\mc{V}| = n$ vertices, $\mc{E} = \{ e_{ij} = (i, j)|~i, j \in \mc{V}, \text{ and } i \neq j \}$ denotes the set of $|\mc{E}| = m$ edges, and $\mc{A} = \{\m{A}_{ij} \in \mb{R}^{d\times d} |~(i,j) \in \mc{E},  \m{A}_{ij}=\m{A}_{ij}^T \geq 0 \}$ denotes the set of matrix weights, one for each edge in $\mc{E}$.\footnote{From the definition, $(i,j)$ and $(j,i)$ denote the same connection between two vertices $i$ and $j$.} The dimension $d$ ($d \geq 1$) of the matrix weights in $\mc{A}$ depends on the problem. Clearly, if $d=1$, the graph $\mc{G}$ becomes an usual undirected scalar-weighted graph. 

Depending on the matrix weights, the interconnection between vertices in $\mc{G}$ are classified into two types. If the matrix weight $\m{A}_{ij}$ corresponding to edge $(i,j) \in \mc{E}$ is positive definite, we say that $(i,j)$ is a positive definite edge and $i$ and $j$ are connected via a positive definite connection. If the weight matrix $\m{A}_{ij}$ corresponding to an edge $(i,j) \in \mc{E}$ is positive semidefinite, we say that $(i,j)$ is a semi-positive definite edge and $i$ and $j$ are connected via a positive definite connection. Apparently, if $i$ and $j$ are disconnected, $\m{A}_{ij} = \m{0}$. We also assume that the interconnections between any two vertices are symmetric, i.e., $\m{A}_{ij} = \m{A}_{ji}$, $\forall (i,j) \in \mc{E}$. 

A path is a sequence of vertices in $\mc{G}$, denoted by $\mc{P} = i_1i_2\ldots i_l$, such that $i_k \neq i_l,$ $\forall i_k, i_l \in \mc{P}$, and each edge $(i_k,i_{k+1})$, $k=1, \ldots, l-1$, is a positive definite or a positive semidefinite connection. The graph $\mc{G}$ is \emph{positive semiconnected} if and only if there exists a path between any two vertices in $\mc{G}$. Otherwise, $\mc{G}$ is disconnected.

In this paper, we mostly focus on positively semiconnected graphs. Graphs with disconnected components can be studied similarly. Assuming that $\mc{G}$ is positive semi-connected, we have the following definitions.

\begin{Definition}[Positive path]\label{def:pos_path} A positive path is a sequence of vertices in $\mc{G}$, denoted by $\mc{P} = i_1i_2\ldots i_l$, such that $i_k \neq i_l,$ $\forall i_k, i_l \in \mc{P}$, and each edge $(i_k,i_{k+1})$, $k=1, \ldots,~l-1$, is a positive definite edge.
\end{Definition}

A tree is an undirected graph containing at least one vertex in which any two vertices are connected by exactly one path. We have the following definition.

\begin{Definition}[Positive tree]\label{def:pos_tree} A positive tree $\mc{T}$ is a tree contained in $\mc{V}$ having all positive connections.
\end{Definition}

Equivalently, for all $i, j \in \mc{T}$, there exists a positive path in $\mc{T}$ connecting $i$ and $j$.

\begin{Definition}[Positive spanning tree]\label{def:pos_span_tree} A positive spanning tree $\mc{T}$ of $\mc{G}$ is a positive tree containing all vertices in $\mc{V}$.
\end{Definition}

\begin{figure}
\begin{center}
\includegraphics[width=7cm]{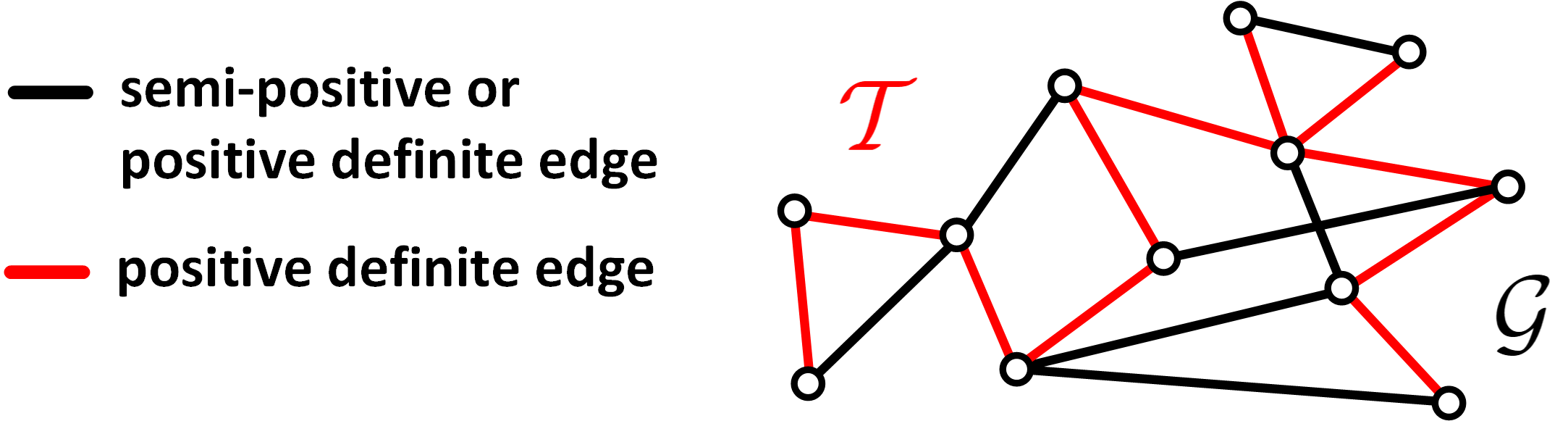}
\caption{$\mc{T}$ is a positive spanning tree of $\mc{G}$. The edges in $\mc{E}(\mc{T})$ are in red color.}
\label{fig:pos-span-tree}
\end{center}
\end{figure}

Note that a tree of $k$ vertices ($k \geq 1$) contains exactly $k-1$ edges. Thus, a positive spanning tree of $\mc{G}$ contains exactly $n-1$ positive connections. An example of positive spanning tree is depicted in Figure \ref{fig:pos-span-tree}. Next, we define several algebraic structures corresponding to the matrix weighted graph $\mc{G}$.
The \emph{matrix-weighted adjacency matrix} of $\mc{G}$ is defined as follows:
\begin{equation}\label{eq:adjacency_matrix}
\m{A} = \begin{bmatrix}
\m{0}&\m{A}_{12}& \cdots &\m{A}_{1n}\\
{{\m{A}_{21}}}&\m{0}& \cdots &\m{A}_{2n}\\
 \vdots & \vdots & \ddots & \vdots \\
{{\m{A}_{n1}}}&{{\m{A}_{n2}}}& \cdots &{{\m{0}}}
\end{bmatrix} \in \mb{R}^{dn \times dn}.
\end{equation}
Since $\mc{G}$ is undirected and $\m{A}_{ij} = \m{A}_{ji}$, it is easy to see that $\m{A}$ is symmetric. For each vertex $i$, the neighbor set of vertex $i$ is defined as $\mc{N}_i = \{ j \in \mc{V}|~(i,j) \in \mc{E} \}$. Let the matrix $\m{D}_i = \sum_{j\in \mc{N}_i} \m{A}_{ij}$ be the \emph{degree matrix} of the vertex $i$. Further, we define $\m{D} = \text{blkdiag}(\m{D}_i)$, the block diagonal matrix of all vertices, as the degree matrix of the graph $\mc{G}$. The \emph{matrix-weighted Laplacian} is defined as follows:
\begin{equation*} 
\m{L} = \m{D} - \m{A} = \begin{bmatrix}
\sum\limits_{j \in \mc{N}_1} \m{A}_{1j}&-\m{A}_{12}& \cdots &-\m{A}_{1n}\\
-{{\m{A}_{21}}}&\sum\limits_{j\in \mc{N}_2} \m{A}_{2j}& \cdots &-\m{A}_{2n}\\
 \vdots & \vdots & \ddots & \vdots \\
-{{\m{A}_{n1}}}&-{{\m{A}_{n2}}}& \cdots &\sum\limits_{j\in \mc{N}_n} \m{A}_{nj}
\end{bmatrix} \in \mb{R}^{dn \times dn}.
\end{equation*}
Consider an arbitrary index of the edges of $\mc{G}$. We can write the edge set and the matrix-weight set as $\mc{E} = \{{e}_{k_{ij}}\}_{k = 1, \ldots, m}$ and $\mc{A} = \{\m{A}_{k_{ij}}\}_{k = 1, \ldots, m}$, correspondingly. From now on, if it is not important to specify the end-vertices explicitly, we will dropout the subscript $ij$ and write $e_k$ and $\m{A}_k$ without ambiguity. 

Let $\m{H} =[h_{ij}] \in \mb{R}^{m \times n}$ denote the incidence matrix corresponding to an arbitrary orientation of the edges in $\mc{E}$. The entries of $\m{H}$ are given as follows:
\[{ {h}_{ki}} = \left\{ {\begin{array}{*{20}{c}}
{\begin{array}{*{20}{c}}
1\\
{ - 1}\\
0
\end{array}}&{\begin{array}{*{20}{l}}
{ \text{if vertex $i$ is the tail of } e_{k},}\\
{ \text{if vertex $i$ is the head of } e_{k},}\\
{\text{otherwise.}}
\end{array}}
\end{array}} \right.\]
An edge $e_k$ is called adjacent to a vertex $i$ if and only if $i$ is a head or a tail of $e_k$, and this adjacency relationship is denoted by $e_k \sim i$.

\begin{Lemma}\label{lem:laplacian_incidence} The matrix-weighted Laplacian can be written in the following form:
\begin{equation} \label{eq:laplacian_incidence}
\m{L} = \bar{\m{H}}^T \rm{blkdiag} (\m{A}_{k}) \bar{\m{H}},
\end{equation}
where $\bar{\m{H}} = (\m{H}\otimes \m{I}_d)$, and $\otimes$ denotes the Kronecker product.
\end{Lemma}

\begin{IEEEproof}
Considering the $ij$-th $d\times d$ block matrix of $\bar{\m{H}}^T \text{blkdiag} (\m{A}_{k_{ij}}) \bar{\m{H}}$, we have 
\[ [\bar{\m{H}}^T \text{blkdiag} (\m{A}_{k}) \bar{\m{H}}]_{ij} = (i\text{-th block column of } \bar{\m{H}})^T \cdot \text{blkdiag} (\m{A}_{k}) \cdot( j\text{-th block column of }\bar{\m{H}}) = \sum_{k=1}^m [\bar{\m{H}}]_{ki}^T \m{A}_{k} [\bar{\m{H}}]_{kj}.\] There are three cases:
\begin{itemize}
\item If $i = j$, since $[\bar{\m{H}}]_{ki} = {h}_{ki} \otimes \m{I}_{d\times d} = {h}_{ki} \m{I}_{d\times d}$, 
\begin{align*}
[\bar{\m{H}}^T \text{diag}(\m{A}_{k}) \bar{\m{H}}]_{ij} = \sum_{k=1}^m [\bar{\m{H}}]_{ki}^T \m{A}_{k} [\bar{\m{H}}]_{ki}= \sum_{k=1}^m (h_{ki} \m{I}_{d\times d}) \m{A}_{k} (h_{ki} \m{I}_{d\times d})= \sum_{k=1}^m (h_{ki})^2 \m{A}_{k} = \sum_{k| e_k \sim i}  \m{A}_{k} = \sum_{j \in \mc{N}_i} \m{A}_{k_{ij}}.
\end{align*}
\item If $i \neq j$ and no edge exists between $i$ and $j$,
\begin{align*}
[\bar{\m{H}}^T \text{diag}(\m{A}_{k}) \bar{\m{H}}]_{ij} &= \sum_{k=1}^m [\bar{\m{H}}]_{ki}^T \m{A}_{k} [\bar{\m{H}}]_{kj}= \sum_{k=1}^m (h_{ki} h_{kj}) \m{A}_{k} = \m{0}.
\end{align*}
\item If $i \neq j$ and $(i,j) \in \mc{E}(\mc{G})$,
\begin{align*}
[\bar{\m{H}}^T \text{diag}(\m{A}_{k}) \bar{\m{H}}]_{ij} &= \sum_{k=1}^m [\bar{\m{H}}]_{ki}^T \m{A}_{k} [\bar{\m{H}}]_{kj}= \sum_{k=1}^m (h_{ki} h_{kj}) \m{A}_{k}= (h_{ki} h_{kj}) \m{A}_{k_{ij}} = -\m{A}_{k_{ij}}.
\end{align*}
\end{itemize}
\end{IEEEproof}
\begin{Corollary}\label{cor:of_lem1} For any vector $\m{v} = [\m{v}_1^T, \ldots, \m{v}_n^T]^T \in \mb{R}^{dn}$, 
$$\m{v}^T\m{L} \m{v} =\sum_{(i,j) \in \mc{E}} (\m{v}_i - \m{v}_j)^T \m{A}_{ij} (\m{v}_i - \m{v}_j).$$
\end{Corollary}
\begin{IEEEproof} We can write
\begin{align}
\m{x}^T\m{L}\m{x} &= \m{x}^T \bar{\m{H}}^T \text{blkdiag}(\m{A}_{k_{ij}}) \bar{\m{H}} \m{x} = (\bar{\m{H}}\m{x})^T \text{blkdiag}(\m{A}_{k_{ij}}) (\bar{\m{H}} \m{x}).
\end{align}
The result follows immediately by observing that the $k$-th block matrix of $\bar{\m{H}}\m{x}$ corresponds to the edge $e_{k}=(i,j) \in \mc{E}$ and is $(\m{x}_j - \m{x}_i)$.
\end{IEEEproof}
\subsection{Matrix-Weighted Consensus Protocols}
Consider a networked dynamic system consisting of $n$ agents. Each agent $i$ in the system has a state vector $\m{x}_i = [x_{i1},\ldots, x_{id}]^T \in \mb{R}^d$, where $d \geq 1$. 
The overall system's states are described by a stacked column vector $\m{x} = [\m{x}_1^T,\ldots,\m{x}_d^T]^T \in \mb{R}^{dn}$. 

The matrix-weighted undirected graph $\mc{G} = (\mc{V}(\mc{G}),\mc{E}(\mc{G}),\mc{A}(\mc{G}))$ describes the interconnection between the agents in the system. Assume that $\mc{G}$ is {positive semi-connected}. An edge $e_{ij} \in \mc{E}$ exists if and only if agent $i$ and agent $j$ can sense their relative state information in at least one state variable. 

In this paper, we consider the agents with single-integrator model. Each agents in the system updates its states under the following protocol:
\begin{equation} \label{eq:consensus_protocol}
\dot{\m{x}}_i = \sum\limits_{j \in \mc{N}_i} \m{A}_{ij} (\m{x}_j  - \m{x}_i), \forall i = 1, \ldots, n,
\end{equation}
where $\m{x}_i \in \mb{R}^d$ is the state and the right-hand side of (\ref{eq:consensus_protocol}) is the control input of agent $i$, $i = 1, \ldots, n$ at time instance $t\geq 0$. Using the matrix-weighted Laplacian, we can express the dynamics of $n$ agents in the following matrix form:
\begin{equation} \label{eq:consensus}
\dot{\m{x}} = -\m{L} \m{x}.
\end{equation}
We have the following definitions.
\begin{Definition}[Consensus] \label{def:consensus} The $n$-agent system is said to  achieve a consensus if and only if $\m{x}_i = \m{x}_j$, for all $i, j \in \mc{V}$, $i \neq j$.
\end{Definition}
Define $\mc{R} = \text{Range} \{ \m{1}_n \otimes \m{I}_{d}\}$ as the consensus space. A consensus of the $n$-agent system is globally/locally asymptotically achieved if and only if $\m{x}(t)$ globally/locally asymptotically approaches $\mc{R}$. Although consensus is important objective, in some applications, the agents' states are desired to converge to some different values. Under the consensus protocol \eqref{eq:consensus_protocol}, clustering behaviors appear naturally. This phenomenon is due to the existence of some positive semidefinite edges in the graph. A partition of $\mc{V}(\mc{G})$ is given by $\mc{C}_1,\ldots, \mc{C}_l,$ $(1\leq l \leq n)$ satisfying two properties: (i) $\mc{C}_i \bigcap \mc{C}_j = \emptyset$, for $i \neq j$, and (ii) $\bigcup_{k=1}^l \mc{C}_k = \mc{V}(\mc{G})$. We have the following definition.

\begin{Definition}[Cluster Consensus] \label{def:stability} The $n$-agent system is said to achieve a clustered consensus if there exists a partition $\mc{C}_1,\ldots, \mc{C}_l,$ such that all agents belonging to the same partition achieve consensus, while for any two agents $i$ and $j$ belonging to two different partitions, $\m{x}_i \neq \m{x}_j$. Each $\mc{C}_i$, $i =1, \ldots, l,$ is referred to as a cluster.
\end{Definition}


\section{Algebraic Condition For Reaching A Consensus}
\label{sec:3}
This section aims to find an algebraic condition of the matrix-weighted Laplacian for reaching consensus. We firstly state several  properties of the matrix-weighted Laplacian and the dynamical system \eqref{eq:consensus}.

\begin{Lemma}[Nullspace of the matrix-weighted Laplacian] \label{lem:laplacian}
The matrix-weighted Laplacian $\m{L}$ is symmetric, positive semidefinite, and $\mc{N}(\m{L}) = \text{span}\{ \m{1}_n \otimes \m{I}_{d}, \{\m{v} = [\m{v}_1^T, \ldots, \m{v}_n^T]^T \in \mb{R}^{dn}| (\m{v}_j - \m{v}_i) \in \mc{N}(\m{A}_{ij}), \forall (i,j) \in \mc{E} \} \}$.
\end{Lemma}
\begin{IEEEproof} The symmetric property of $\m{L}$ follows immediately from its definition. From \eqref{eq:laplacian_incidence}, we can write
\begin{align} \label{eq:m(t)}
\m{L} &= \bar{\m{H}}^T \text{blkdiag} (\m{A}_k) \bar{\m{H}} = \bar{\m{H}}^T \text{blkdiag} (\m{A}_k^{1/2}) \text{ blkdiag} (\m{A}_k^{1/2}) \bar{\m{H}}  = \m{M}^T\m{M}, 
\end{align}
where $\m{M} = \text{ blkdiag} (\m{A}_k^{1/2}) \bar{\m{H}}$. Equation \eqref{eq:m(t)} shows that $\m{L}$ is positive semidefinite. Moreover, we have $\mc{N}(\m{L}) = \mc{N}(\m{M}^T\m{M})=\mc{N}(\m{M})$. As a result, $\mc{N}(\m{L}) = \mathcal{N}{(\m{M})} \supseteq \mathcal{N}(\bar{\m{H}}) = \text{span} \{ \m{1}_n \otimes \m{I}_{d} \}$. Consider $\m{v} = [\m{v}_1^T, \ldots, \m{v}_n^T]^T \notin \text{ span }\{ \m{1}_n \otimes \m{I}_{d}\}$ such that $\m{L} \m{v} = \m{0} $. It follows $\m{v}^T\m{L} \m{v} = 0$. Thus, from Corollary \ref{cor:of_lem1}, we have 
\begin{equation} \label{eq:null-2}
\sum_{(i,j) \in \mc{E}} (\m{v}_i - \m{v}_j)^T \m{A}_{ij} (\m{v}_i - \m{v}_j) = 0.
\end{equation}
Since $\m{A}_{ij}$ is symmetric and positive semidefinite, \eqref{eq:null-2} implies that $(\m{v}_i - \m{v}_j) \in \mc{N}(\m{A}_{ij})$, for all $(i,j) \in \mc{E}$. This concludes the proof.
\end{IEEEproof}

\begin{Remark}\label{remark:zero_eig_of_L}
Based on Lemma \ref{lem:laplacian}, it follows that $\rm{dim}( \mc{N}(\m{L}) ) \geq \rm{dim} (\mc{R})$. Thus, the matrix weighted Laplacian $\m{L}$ has at least $d$ zero eigenvalues. Let $\{\lambda_i\}_{i=1,\ldots, dn}$ be the eigenspectra of $\m{L}$; then we have $0 = \lambda_1 = \ldots = \lambda_d \leq \lambda_{d+1} \leq \ldots \leq \lambda_{dn}$.
\end{Remark}

\begin{Lemma}\label{lem:centroid_invariance}
Under the consensus protocol \eqref{eq:consensus_protocol}, the average $\bar{\m{x}} = \frac{1}{n}\sum_{i=1}^n \m{x}_{i}$ is invariant.
\end{Lemma}
\begin{IEEEproof}
The average state can be written as $\bar{\m{x}} = \frac{1}{n}(\m{1}_n^T \otimes \m{I}_{d}) \m{x}$. Taking the derivative of $\bar{\m{x}}$ along the trajectory of \eqref{eq:consensus} yields
\begin{align}
\dot{\bar{\m{x}}} &= \frac{1}{n}(\m{1}_n^T \otimes \m{I}_{d}) \dot{\m{x}} = -\frac{1}{n}(\m{1}_n^T \otimes \m{I}_{d}) {\m{L}} \m{x}.
\end{align}
Since $\m{L}$ is symmetric, if $\m{v} \in \mc{N}(\m{L})$, then $\m{v}^T$ belongs to the right nullspace of $\m{L}$. As a result, $(\m{1}_n^T \otimes \m{I}_{d})  {\m{L}} = \m{0}$ and it follows $\dot{\bar{\m{x}}} = \m{0}$, i.e. the system's average is invariant.
\end{IEEEproof}

The following theorem characterizes the dynamical behavior of the consensus protocol \eqref{eq:consensus}. 
\begin{Theorem}[Stability]\label{thm:stability} Assume that $\mc{G}$ is positive semi-connected. Then any trajectory of \eqref{eq:consensus} asymptotically approaches the invariant set $\mc{N}(\m{L})$.
\end{Theorem}
\begin{IEEEproof}
Consider the potential function $V = \frac{1}{2} \|\m{x}\|^2$, 
which is positive definite, radially unbounded, and continuously differentiable. The derivative of $V$ along the trajectory of \eqref{eq:consensus} is given by
\begin{align*}
\dot{V} &= -\m{x}^T \m{L} \m{x} = -\sum_{(i,j) \in \mc{E}(\mc{G})} (\m{x}_i - \m{x}_j)^T \m{A}_{ij} (\m{x}_i - \m{x}_j) \leq 0.
\end{align*}
It follows that $\|\m{x}\| \leq \|\m{x}(0)\|$, or i.e., $\m{x}(t)$ is bounded. Further, $\dot{V}$ is negative semidefinite and $\dot{V}=0$ if and only if $\m{x} \in \mc{N}(\m{L})$. Based on LaSalle's invariance principle, any trajectory of \eqref{eq:consensus} asymptotically approaches the invariant set $\mc{N}(\m{L})$ as described in Lemma \ref{lem:laplacian}.
\end{IEEEproof}

\begin{Lemma} \label{lem:unique_eq_pt} If $\mc{N}(\m{L}) = \rm{span}\{\m{1}_n \otimes \m{I}_{d}\}$, the system \eqref{eq:consensus} has a unique equilibrium point $\m{x}^* = \m{1}_n \otimes \bar{\m{x}}$.
\end{Lemma}
\begin{IEEEproof} We prove that $\m{1}_n \otimes \bar{\m{x}}$ is the unique equilibrium of \eqref{eq:consensus} if $\mc{N}(\m{L})= span\{\m{1}_n \otimes \m{I}_{d}\}$ by contradiction. Let $\{\m{e}_i\}_{i =1, \ldots, d}$ be a basis of $\mb{R}^d$, where $\m{e}_i = [0, \ldots, 1, \ldots, 0]^T$ is a vector with all zero entries except for an $1$ on the $i$th row. Suppose that there exists $\m{x}' \in \text{span} \{ \m{1}_n \otimes \m{I}_{d} \}$ such that $\m{x}' \neq \m{x}^*$. Since $\m{x}' \in \text{span} \{ \m{1}_n \otimes \m{I}_{d} \}$, we write
$$\m{x}' = \sum_{i=1}^d \bar{x}'_i (\m{1}_n \otimes \m{e}_i) = \m{1}_n \otimes \bar{\m{x}}' = (\m{1}_n \otimes \m{I}_{d})\bar{\m{x}}',$$
where $\bar{\m{x}}' = [\bar{x}_1', \ldots, \bar{x}_d']^T$. It follows from Lemma \ref{lem:centroid_invariance} that
\begin{align*}
\bar{\m{x}} = \frac{1}{n} (\m{1}_n^T \otimes \m{I}_{d}) \m{x}'   = \frac{1}{n} (\m{1}_n^T \otimes \m{I}_{d}) (\m{1}_n \otimes \m{I}_{d})\m{\bar{x}}'= \frac{1}{n} (\m{1}_n^T\m{1}_n \otimes \m{I}_{d}) \m{\bar{x}}' = \frac{1}{n} (n \otimes \m{I}_{d}) \m{\bar{x}}'  = \m{\bar{x}}'.
\end{align*}
Thus, $\m{x}' = \m{1}_n \otimes \bar{\m{x}} = \m{x}^*$, which is a contradiction. This contradiction implies that $\m{x}^* = \m{1}_n \otimes \bar{\m{x}}$ is the unique equilibrium of \eqref{eq:consensus}.
\end{IEEEproof}

The following theorem gives a necessary and sufficient condition for \eqref{eq:consensus} to globally achieve an average consensus.

\begin{Theorem}[Average Consensus] \label{thm:consensus_condition} The system \eqref{eq:consensus} globally exponentially converges to the system's average  $\m{x}^*=\m{1}_n\otimes \bar{\m{x}}$ if and only if $\mc{N}(\m{L}) = \text{span}\{\m{1}_n \otimes \m{I}_{d}\}$.
\end{Theorem}

\begin{IEEEproof}(Necessity). We prove by contradiction. Assume that \eqref{eq:consensus} globally asymptotically converges to $\m{x}^*=\m{1}_n \otimes \bar{\m{x}}$ but $\mc{N}(\m{L}) \neq \mc{R}$. From Lemma \ref{lem:laplacian}, there exists $\m{x}' \in \mb{R}^{dn}$ such that $\m{L}\m{x}'=\m{0}$ and $\m{x}' \notin \mc{R}$. Thus, $\m{x} = \m{x}'$ is also an equilibrium point of \eqref{eq:consensus}, and any trajectory with $\m{x}(0) = \m{x}'$ stays at $\m{x}'$ for all $t \geq 0$. Thus, $\m{x}^*$ is not globally asymptotically stable, which contradicts the assumption.

(Sufficiency). Suppose that $\mc{N}(\m{L}) = \mc{R}$. Following the proof of Theorem \ref{thm:stability}, any trajectory of \eqref{eq:consensus} converges to $\mc{N}(\m{L}) = \{ \m{1}_n \otimes \m{I}_{d} \}$. It follows from Lemma \ref{lem:unique_eq_pt} that $\m{x}^* = \m{1}_n \otimes \bar{\m{x}} \in \mc{N}(\m{L})$ is the unique equilibrium point of \eqref{eq:consensus}. 

Consider the potential function $V = \frac{1}{2} \bm{\delta}^T\bm{\delta}$, where $\bm{\delta} = \m{x} - \m{1}_n \otimes \bar{\m{x}}$ is the disagreement vector. Then, V is positive definite, radially unbounded, and continuously differentiable. The derivative of $V$ along the trajectory of \eqref{eq:consensus} is
\begin{align} \label{eq:dot_V1}
\dot{V} =\bm{\delta}^T \dot{\bm{\delta}} = - \bm{\delta}^T\m{L} \m{x} = -\bm{\delta}^T \m{L} \bm{\delta} \leq 0,
\end{align}
where in the third equality, we have used the fact that
$\m{L} \bm{\delta} = \m{L} \m{x} - \m{L} (\m{1}_n \otimes \bar{\m{x}}) = \m{L} \m{x} - \m{L} (\m{1}_n \otimes \m{I}_{d}) \bar{\m{x}} = \m{L} \m{x}$. 
%

Moreover, $\bm{\delta} \perp \mc{R}$ since $(\m{1}_n \otimes \m{I}_{d})^T \bm{\delta} = (\m{1}_n \otimes \m{I}_{d})^T \m{x} - (\m{1}_n^T\m{1}_n \otimes \m{I}_{d}) \bar{\m{x}} = n \bar{\m{x}} - n \bar{\m{x}} = \m{0}$. Therefore, we can write 
\begin{align}\label{eq:convergence_rate}
\dot{V} = -\bm{\delta}^T \m{L} \bm{\delta}  \leq - \lambda_{d+1}(\m{L}) \bm{\delta}^T\bm{\delta} \leq -\alpha V \leq 0,
\end{align}
where $\alpha = 2 \lambda_{d+1}(\m{L}) > 0$. Further, $\dot{V} = 0$ if and only if $\bm{\delta} = \m{0}$, or $\m{x} = \m{x}^*= \m{1}_n \otimes \bar{\m{x}}$. Therefore, the equilibrium $\m{x}^*$ is globally exponentially stable, i.e. \eqref{eq:consensus} globally exponentially achieves an average consensus.
\end{IEEEproof}

\begin{Remark}\label{remark:convergence_rate}
Equation \eqref{eq:convergence_rate} shows that $\lambda_{d+1}$, the smallest positive eigenvalue of $\m{L}$, determines the convergence rate of the matrix-weighted consensus protocol \eqref{eq:consensus}. Thus, $\lambda_{d+1}$ is a performance index of the network, and this index is analogous to the algebraic connectivity of $\mc{G}$ in the usual consensus algorithm \cite{Olfati2007,MesbahiEgerstedt}.
\end{Remark}

In the usual consensus algorithm, the average consensus is asymptotically achieved if and only if the graph is connected and the weights are positive scalars \cite{Olfati2004,Olfati2007}. Thus, we expect \eqref{eq:consensus} to reach a consensus when all matrix weights are positive definite.
\begin{Corollary} \label{cor:of_thm_2} Under the consensus protocol \eqref{eq:consensus_protocol}, if $\m{A}_{ij} > 0$ for all $(i,j) \in \mc{E}$, all agents globally exponentially achieve a consensus.
\end{Corollary}

\begin{IEEEproof}  
Since $\m{A}_{k}$, $k = 1, \ldots, m,$ are positive definite from the assumption, it follows that $\mc{N}(\m{L}) = \mc{N}(\m{M}) = \mc{N}(\text{diag}(\m{A}_{k}^{1/2})\bar{\m{H}})= \mc{N}(\bar{\m{H}}) = \{ \m{1}_n \otimes \m{I}_{d} \}$. Thus, it follows from Theorem \ref{thm:consensus_condition} that \eqref{eq:consensus} globally exponentially achieves a consensus. 
\end{IEEEproof}

\section{Algebraic Graph Theory of Consensus and Clustered Consensus Phenomena}
\label{sec:4}
In the previous section, Theorem \ref{thm:consensus_condition} provides an algebraic condition for reaching a consensus. However, that condition requires finding the nullspace of $\m{L}$. Corollary \ref{cor:of_thm_2} gives a sufficient condition for achieving consensus. The condition is a quite clear and straightforward. However, since the condition is only sufficient, it might be conservative. In this section, we aim to find some conditions for consensus and clustered consensus related with the  matrix-weighted graph $\mc{G}$.

\begin{Lemma} \label{lem:pos_span_tree}
If there exists a positive spanning tree in $\mc{G}$, then an average consensus is globally exponentially achieved.
\end{Lemma}

\begin{IEEEproof}
Suppose $\mc{G}$ has a spanning tree $\mc{T}$ having all edges with positive definite matrix weights. We can label the edges of $\mc{G}$ such that the $n-1$ edges in $\mc{T}$ are $e_1, e_2, \ldots, e_{n-1}$ and the remaining $m-n+1$ edges in $\mc{E}$ are $e_{n}, e_{n+1}, \ldots, e_m$. The incidence matrix corresponding to this labeling can be written as
\begin{equation*}
\m{H} = \begin{bmatrix}
\m{H}_{\mc{E}(\mc{T})} \\
\m{H}_{\mc{E}\setminus \mc{E}(\mc{T})}
\end{bmatrix},
\end{equation*}
where $\m{H}_{\mc{E}(\mc{T})} \in \mb{R}^{(n-1)\times n}$ represents $n-1$ edges of $\mc{T}$ and $\m{H}_{\mc{E}\setminus \mc{E}(\mc{T})}
\in \mb{R}^{(m-n+1)\times n}$ represents the remaining edges in the graph. Note that the rows of $\m{H}_{\mc{E}\setminus \mc{E}(\mc{T})}
$ are linearly dependent on the rows of $\m{H}_{\mc{E}(\mc{T})}$ \cite{Zelazo2011tac}. Specifically, there exists a matrix $\m{T} \in \mb{R}^{(m-n+1) \times m}$ such that:
$$ \m{T} \m{H}_{\mc{E}(\mc{T})} =\m{H}_{\mc{E}\setminus \mc{E}(\mc{T})},$$
where 
$ \m{T} = \m{H}_{\mc{E}\setminus \mc{E}(\mc{T})} \m{H}_{\mc{E}(\mc{T})}^T (\m{H}_{\mc{E}(\mc{T})}\m{H}_{\mc{E}(\mc{T})}^T)^{-1}.$ 
Thus, we can rewrite the incidence matrix as
\begin{equation}
\m{H} = \begin{bmatrix}
\m{H}_{\mc{E}(\mc{T})} \\
\m{T} \m{H}_{\mc{E}(\mc{T})}
\end{bmatrix} = \begin{bmatrix}
\m{I}_{n-1} \\
\m{T} 
\end{bmatrix} \m{H}_{\mc{E}(\mc{T})},
\end{equation}
Any equilibrium point of \eqref{eq:consensus} must satisfy
\begin{align}
\dot{\m{x}} = - \bar{\m{H}}^T \text{ blkdiag} (\m{A}_k) \bar{\m{H}} \m{x} = \m{0}.
\end{align}
It follows that
$\m{x}^T\bar{\m{H}}^T \text{blkdiag} (\m{A}_k) \bar{\m{H}} \m{x} =0$, or $\|\text{blkdiag} (\m{A}_k^{{1}/{2}}) \bar{\m{H}} \m{x}\|^2 =\|\m{M} \m{x}\|^2 = 0$. Denoting $\bar{\m{T}} = \m{T} \otimes \m{I}_d$, this equation is equivalent to
\begin{equation} \label{eq:tree-graph}
\m{M}\m{x} = \begin{bmatrix}
\text{ blkdiag}(\m{A}_k^{{1}/{2}})_{k=1}^{n-1} \bar{\m{H}}_{\mc{E}(\mc{T})} \m{x} \\
\text{ blkdiag}(\m{A}_k^{{1}/{2}})_{k=n}^{m} \bar{\m{T}} \bar{\m{H}}_{\mc{E}(\mc{T})} \m{x}
\end{bmatrix} = \m{0}.
\end{equation}
Observe that $\text{ blkdiag}(\m{A}_k^{{1}/{2}})\bar{\m{H}}_{\mc{E}(\mc{T})}\m{x} = \m{0}$ is equivalent to $\bar{\m{H}}_{\mc{E}(\mc{T})}\m{x} = \m{0}$ since $\m{A}_k$, $k=1, \ldots, n-1$, are positive definite (the corresponding edges are in the positive spanning tree). Further, since $\m{H}_{\mc{E}(\mc{T})}$ is the incidence matrix corresponding to a tree, we have $\mc{N}(\m{H}_{\mc{E}(\mc{T})}) = \text{span} \{ \m{1}_n \}$, which means $\mc{N}(\m{\bar{H}}_{\mc{E}(\mc{T})}) = \text{span}\{ \m{1}_n \otimes \m{I}_{d}\} = \mc{R}$. It follows from Lemma \ref{lem:unique_eq_pt} that the equilibrium is unique and is $\m{x}^* = \m{1}_n \otimes \bar{\m{x}}$. 
Also, it is easy to check that $\text{blkdiag}(\m{A}_k^{{1}/{2}})_{k=n}^{m} \bar{\m{T}} \bar{\m{H}}_{\mc{E}(\mc{T})} \m{x}^* = \m{0}$.

Finally, the stability of $\m{x}=\m{x}^*$ follows from Theorem \ref{thm:consensus_condition}.
\end{IEEEproof}

\begin{Lemma} \label{lem:pos_tree}
Suppose there exists a positive tree $\mc{T} \subset \mc{G}$ of $l$ vertices. Under the consensus protocol \ref{eq:consensus_protocol}, $\m{x}_i(t) \to \m{x}_j(t)$, $\forall i, j \in \mc{T}$, as $t \to \infty$.
\end{Lemma}

\begin{IEEEproof}
Let the state vector be indexed as $\m{x} = [\m{x}_{\mc{T}}^T, \m{x}_{\mc{V}\setminus \mc{V}(\mc{T})}^T]^T$. We express the incidence matrix in the following form
\[\begin{array}{*{20}{c}}
{}&\begin{array}{*{20}{c}}
{\mc{V}(\mc{T})}&{\rm{   }}\mc{V}(\mc{G}) \setminus \mc{V}(\mc{T}) 
\end{array}\\
{\begin{array}{*{20}{c}}
{{\cal E}\left( {\cal T} \right)}\\
{{\cal E}\left( {{\cal V}\left( {\cal T} \right)} \right)\backslash {\cal E}\left( {\cal T} \right)}\\
{{\cal E}\backslash {\cal E}\left( {{\cal V}\left( {\cal T} \right)} \right)}
\end{array}}&{\left[ {\begin{array}{*{20}{c}}
{{\rm{     }}{\m{H}_1}{\rm{    }}}& {{\rm{        }}\m{0}{\rm{        }}}\\
{{\rm{  }}{\m{H}_2}}& {{\rm{ }} \m{0} {\rm{ }}}\\
{{\rm{  }}{\m{H}_3}}& {{\rm{     }}{\m{H}_4}{\rm{  }}}
\end{array}} \right]}
\end{array}\begin{array}{*{20}{c}}
{}\\
{}\\
{ = \m{H}},\\
{}
\end{array}\]
where $[\m{H}_{1} \quad \m{0}] \in \mb{R}^{(l-1)\times n}$ associates with the $l$ edges belonging to the tree $\mc{T}$, $[\m{H}_{2} \quad \m{0}]$ associates with the $r-l+1 $ edges between vertices in $\mc{V}(\mc{T})$ which do not belong to the tree, and $[\m{H}_3 \quad \m{H}_4]$ associates with the remaining edges in $\mc{E}$. Similarly to the proof of Lemma \ref{lem:pos_span_tree}, $\m{H}_2$ is linearly dependent on $\m{H}_1$ and this dependency is characterized by $\m{H}_2 = \m{T} \m{H}_1$. Therefore, the equilibrium set of \eqref{eq:consensus} must satisfy
\begin{equation}
\m{M}\m{x} = \begin{bmatrix}
\text{ blkdiag}(\m{A}_k^{{1}/{2}})_{k=1}^{l-1} \bar{\m{H}}_{1} \m{x}_{\mc{T}} \\
\text{ blkdiag}(\m{A}_k^{{1}/{2}})_{k=l}^{r} \bar{\m{T}} \bar{\m{H}}_{1} \m{x}_{\mc{T}} \\
\text{ blkdiag}(\m{A}_k^{{1}/{2}})_{k=r+1}^{m} (\bar{\m{H}}_3 \m{x}_{\mc{T}} + \bar{\m{H}}_4 \m{x}_{\mc{V}\setminus \mc{V}(\mc{T})})
\end{bmatrix} = \m{0}.
\end{equation}

Since $\mc{T}$ is a positive tree, $\text{ blkdiag}(\m{A}_k^{{1}/{2}})_{k=1}^{l-1}$ is positive definite. It follows $\bar{\m{H}}_1\m{x} = \m{0}$. Further, since $\m{H}_1$ is the incidence matrix associated with a tree, $\mc{N}(\bar{\m{H}}_1) = \{ \m{1}_{l} \otimes \m{I}_{d} \}$. Hence, any equilibrium $\m{x}^*$ of \eqref{eq:consensus} must have $\m{x}_{\mc{T}}^* \in \mc{N}(\bar{\m{H}}_1)$, i.e., all equilibrium states of $l$ agents belonging to the positive tree $\mc{T}$ are the same. Based on Theorem \ref{thm:stability}, the agents in $\mc{T}$  asymptotically reach a consensus.
\end{IEEEproof}

The following result provides a condition to determine whether or not two vertices belong to a same cluster.
\begin{Theorem} \label{prop:cluster}
Given a positive tree $\mc{T}$, let the cluster $\mc{C}(\mc{T})$ generated from $\mc{T}$ be containing:
\begin{itemize}
\item[i.] all vertices in $\mc{T}$,
\item[ii.] any vertex $i \notin \mc{T}$, which defines the set $\mc{S}_i = \{ \mc{P}_k=\{v^k_1\ldots v^k_{|\mc{P}_k|}\} | v^k_1 = i, v^k_{|\mc{P}_k|}  \in \mc{T}, \text{ and } \forall j = 1,\ldots, {|\mc{P}_k|-1}, v^k_j \notin \mc{T} \}$, satisfying the following conditions:
\begin{itemize}
\item[a.] for each path $\mc{P}_k$, denoting $\mc{N}(\mc{P}_k) = \bigcup_{j=1}^{{|\mc{P}_k|-1}} \mc{N}( \m{A}_{v^k_jv^k_{j+1}})$, it holds
\begin{equation} \label{eq:same_cluster}
\rm{dim} (\bigcap\nolimits_{k=1}^{|\mc{S}_i|} \mc{N}(\mc{P}_k)) = 0
\end{equation}
\item[b.] each path $\mc{P}_{k} \in \mc{S}_i$ has no loop, i.e. $v_l \neq v_m,$ $\forall v_l, v_m \in \mc{P}_k$
\end{itemize}
\end{itemize}
Then, under the consensus protocol \eqref{eq:consensus_protocol}, all agents in the cluster $\mc{C}(\mc{T})$ have the same equilibrium state. Furthermore, in algorithmic perspective, the set $\mc{S}_i$ is finite.
\end{Theorem}
\begin{IEEEproof}
First, all vertices in $\mc{T}$ converge to a same value due to Lemma \ref{lem:pos_tree}. Denote that the common value by $\m{x}^*_{\mc{T}}$. 

Next, consider a vertex $i \notin \mc{T}$ satisfying the condition (ii). Now, we consider the condition (ii.a). Let $\m{x}^*_i$ be the equilibrium state of agent $i$. Then, from the definition of $\mc{N}(\mc{P}_k)$, we can write 
\begin{equation}\label{eq:same_cluster1}
\m{x}_i^* - \m{x}_{\mc{T}}^* \in \mc{N}(\mc{P}_k), \quad \forall \mc{P}_k \in \mc{S}_i.
\end{equation}
It follows from \eqref{eq:same_cluster} that the only solution for $|\mc{S}_i|$ equations \eqref{eq:same_cluster1} is $\m{x}_i^* - \m{x}_{\mc{T}}^* = \m{0}$, or $\m{x}_i^* = \m{x}_{\mc{T}}^*$. Thus, the cluster $\mc{C}(\mc{T})$ reaches a consensus. But, in the above condition (ii.a), there could be infinitely many paths from $i$ to $\mc{T}$ if there are loops. The condition (ii.b) ensures that the number of paths from $i$ to $\mc{T}$ is finite. To show this, suppose that $\mc{P}_1$ and $\mc{P}_2$ are two paths, and $\mc{P}_2$ is obtained by adding loops to $\mc{P}_1$. Then, $\mc{N}(\mc{P}_1) \subseteq \mc{N}(\mc{P}_2)$. It follows $\mc{N}(\mc{P}_1) \cap \mc{N}(\mc{P}_2) = \mc{N}(\mc{P}_1)$; thus it is not necessary to consider loops when checking the condition \eqref{eq:same_cluster}, which means that $\vert \mc{S}_i\vert$ is finite. 

Finally, consider a vertex $i$ which does not satisfy both (i) and (ii). Let $\mc{S}_i$ be the set of paths from $i$ to $\mc{T}$ with $\rm{dim} (\bigcap\nolimits_{k=1}^{|\mc{S}_i|} \mc{N}(\mc{P}_k)) \geq 1$. Then clearly there exists a nontrivial solution, i.e., $\m{x}_i^* - \m{x}_{\mc{T}}^* \neq 0$.   
\end{IEEEproof}

The following result follows immediately from Proposition \ref{prop:cluster}.
\begin{Corollary} Consider a positive tree $\mc{T}$ and a vertex $i \notin \mc{T}$. From $i$ to $j \in \mc{T}$, if there exist at least two paths, which include positive semi-definite weighting matrices, such that eq. (15) holds, then agent $i$ can be added into the cluster $\mc{C}(\mc{T})$.
\end{Corollary}

\begin{Example}\label{ex:1} To illustrate Proposition \ref{prop:cluster}, consider a four-agent system in $\mb{R}^{3}$ with the interaction graph as depicted in Figure~\ref{fig:ex1}. The matrix-weight corresponding to each connections between the agents in the system are given by $\m{A}_{12}=\begin{bmatrix} 0&0&0\\0&1&0\\0&0&1 \end{bmatrix}, $
$\m{A}_{13}=\begin{bmatrix} 1&0&0\\0&0&0\\0&0&0 \end{bmatrix}, $ and $\m{A}_{23}=\begin{bmatrix} 1&0&0\\0&0&0\\0&0&1 \end{bmatrix}$ and $\m{A}_{14} = \begin{bmatrix} 1&0&0\\0&2&0\\0&0&1 \end{bmatrix}$. It is easy to see that $\m{A}_{14}$ is positive definite while other matrix weights are positive semidefinite. As a result, there is a positive tree $\mc{T}$ in the graph containing vertex 1 and vertex 4. Moreover, we have $\mc{N}(\m{A}_{12}) = \text{span} \{ [1,0,0]^T \}$, $\mc{N}(\m{A}_{13}) = \text{span} \{ [0,1,0]^T, [0,0,1]^T\}$, and $\mc{N}(\m{A}_{23}) = \text{span} \{ [0,1,0]^T \}$.

There are two paths (without loop) from vertex 2 to vertex 1 (also to the tree $\mc{T}$): $\mc{P}_1 = 21$, and $\mc{P}_2 = 231$. By definition, we have $\mc{N}(\mc{P}_1) = \mc{N}(\m{A}_{12}) = \text{span} \{ [1,0,0]^T \}$, and $\mc{N}(\mc{P}_2) =  \mc{N}(\m{A}_{13}) \bigcup \mc{N}(\m{A}_{23}) = \text{span} \{ [0,1,0]^T, [0,0,1]^T\}$. It follows $\mc{N}(\mc{P}_1) \bigcap \mc{N}(\mc{P}_2) = \{ \m{0} \}$, which further implies that agent 2 is in the same cluster  $\mc{C}(\mc{T})$ due to Proposition \ref{prop:cluster} (ii).

On the other hand, consider the vertex 3. There are two paths from vertex 3 to the cluster $\mc{C} = \{1, 2, 4\}$: $\mc{P}_3 = 31$ and $\mc{P}_4 = 32$. Since $\mc{N}(\mc{P}_3) \bigcap \mc{N}(\mc{P}_4) = \mc{N}(\m{A}_{23}) = \text{span} \{ [0,1,0]^T \}$, the vertex 3 does not belong to the cluster $\mc{C}$.
\begin{figure}[t]
\centering
\includegraphics[width=6.8cm]{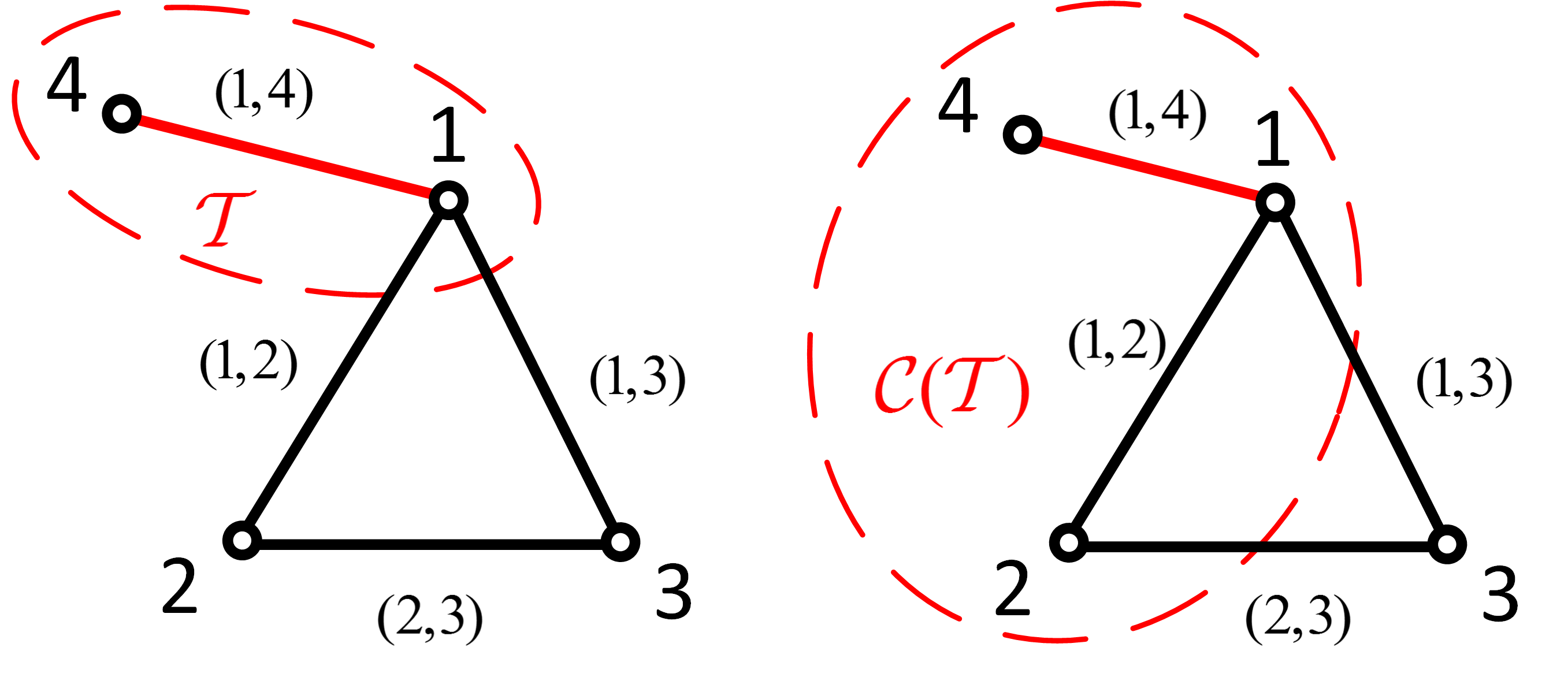}
\caption{Illustration of the four-agent system in Example \ref{ex:1}.}
\label{fig:ex1}
\end{figure}
\begin{figure*}[ht!]
    \begin{subfigure}[b]{0.33\textwidth}
    \centering
    \includegraphics[height=4.5cm]{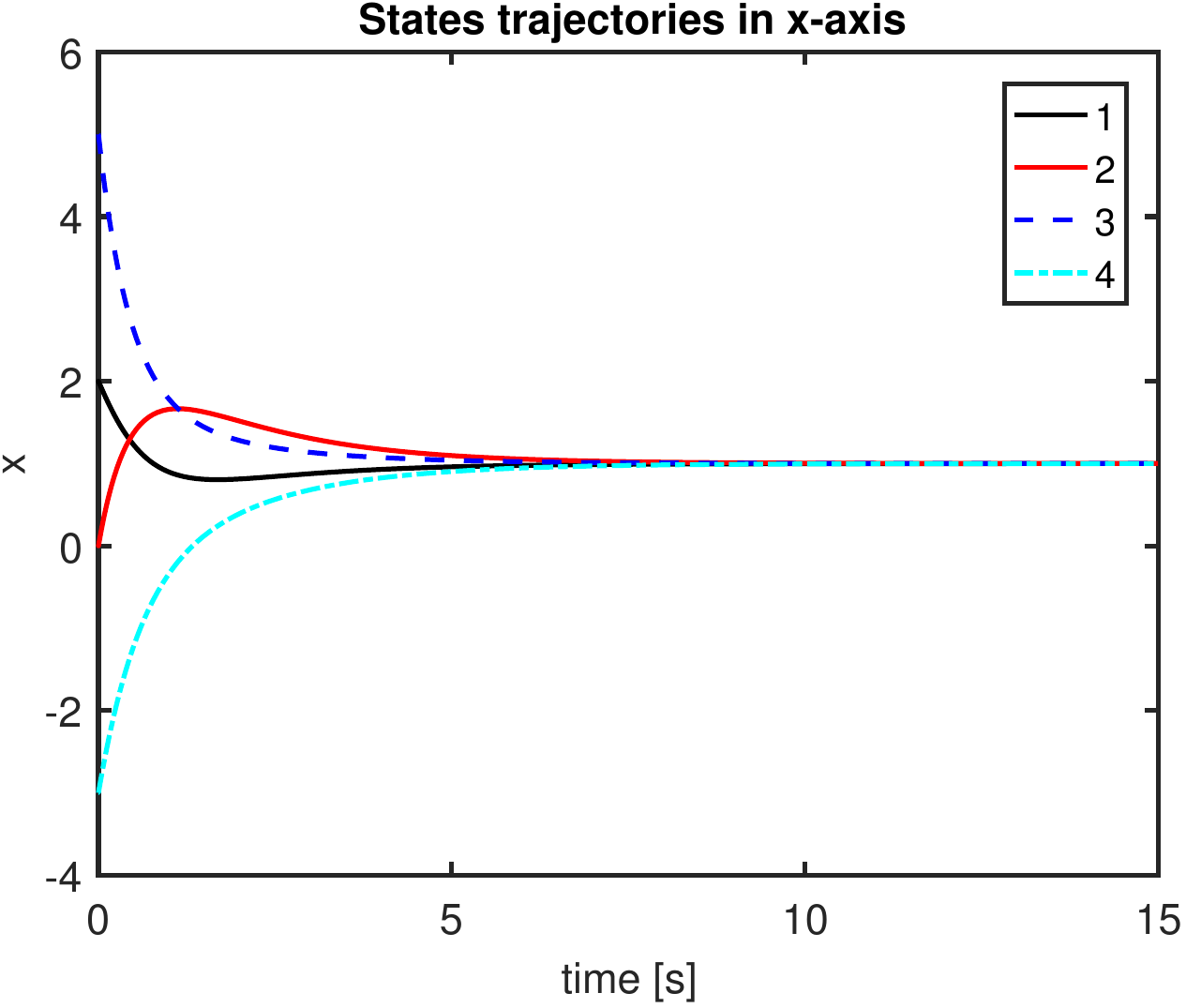}
    \caption{The $x$-axis dynamics.}
    \label{fig:sim-exp1}
    \end{subfigure}
    \begin{subfigure}[b]{0.33\textwidth}
    \centering
    \includegraphics[height=4.5cm]{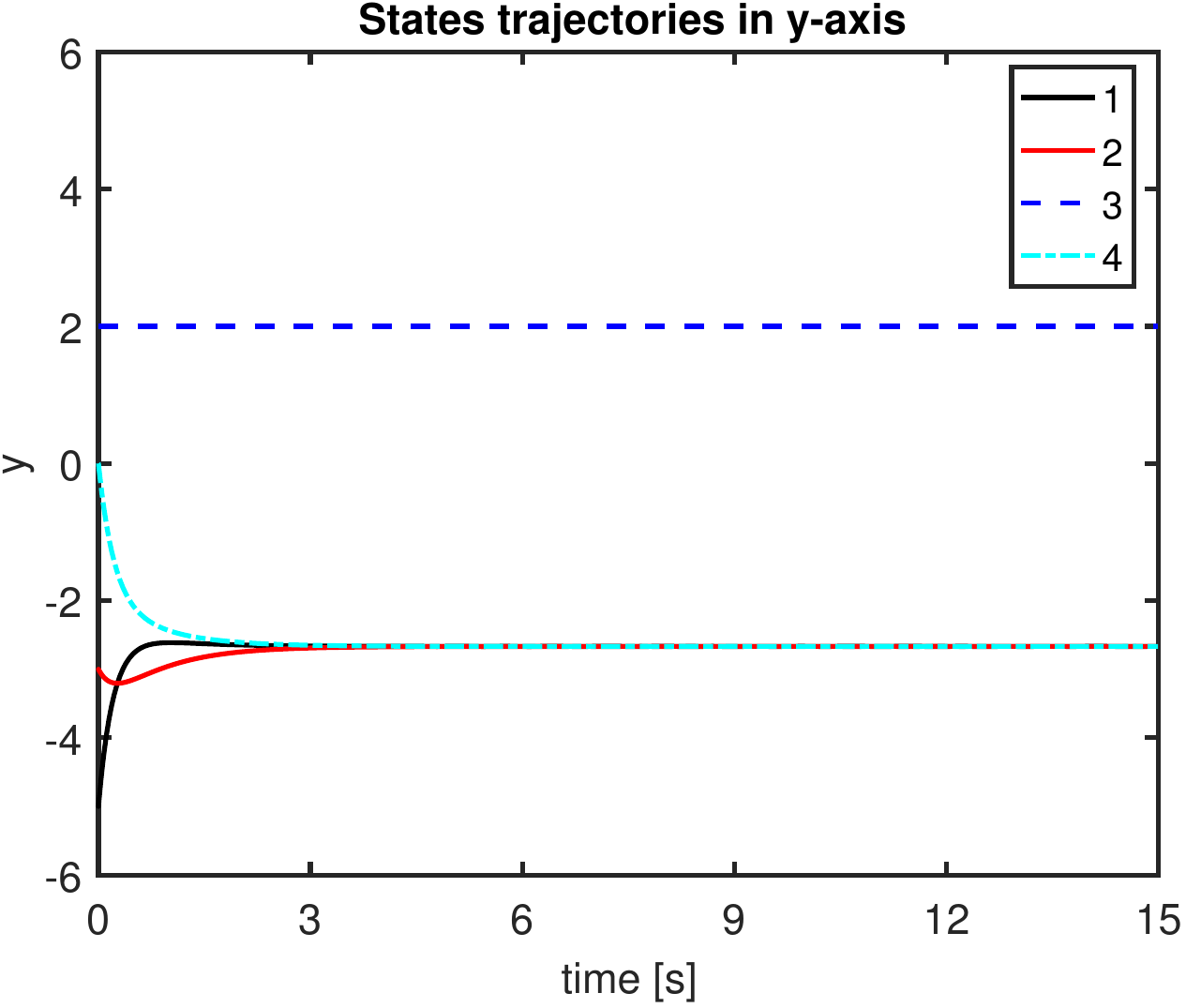}
    \caption{The $y$-axis dynamics.}
    \label{fig:sim-exp2}
    \end{subfigure}    
    \begin{subfigure}[b]{0.33\textwidth}
    \centering
    \includegraphics[height=4.5cm]{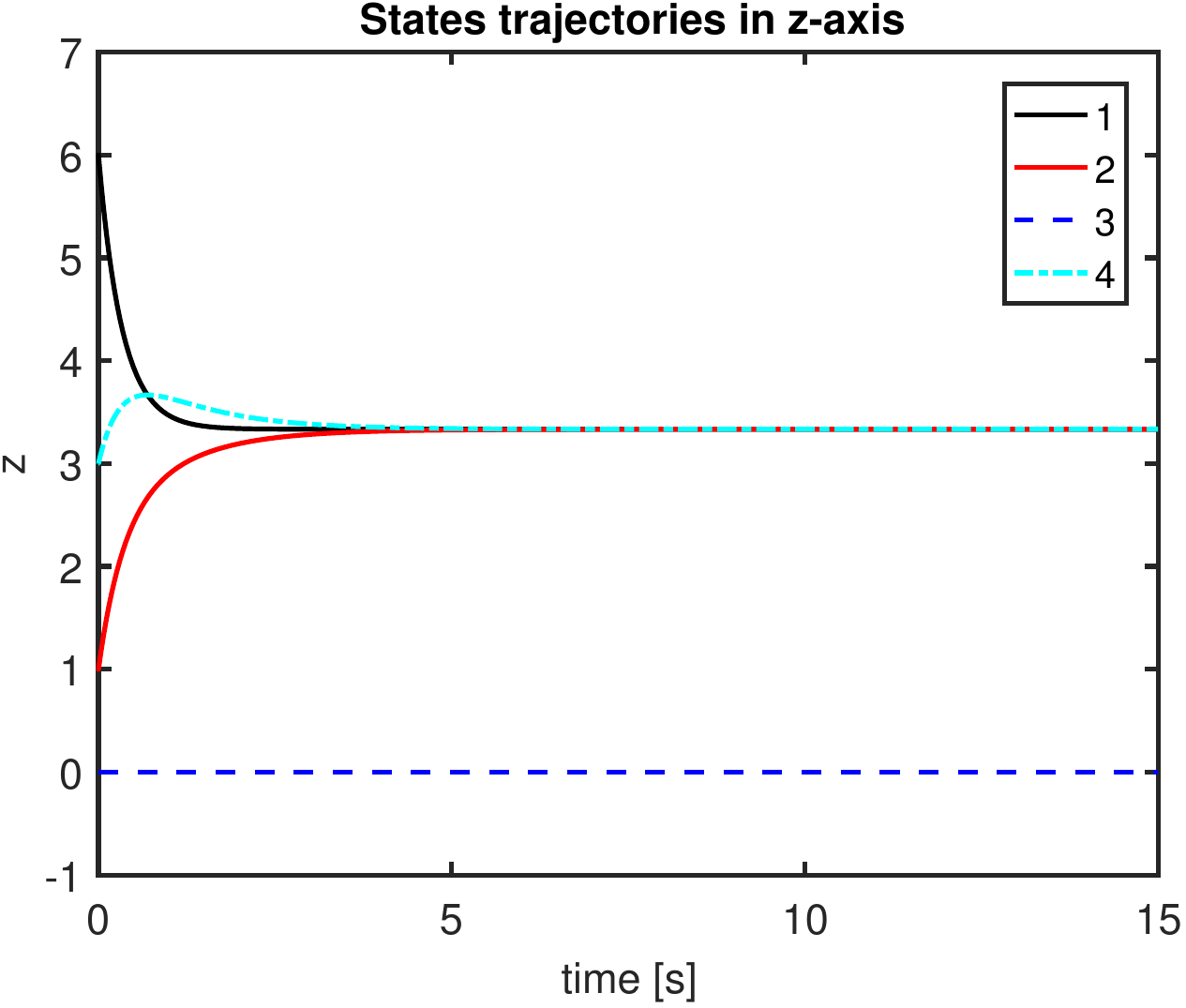}
    \caption{The $z$-axis dynamics.}
    \label{fig:sim-exp3}
    \end{subfigure}
    \caption{Simulation: The states' dynamics of four agents under the consensus protocol \eqref{eq:consensus_protocol}.}
    \label{fig:sim_ex}
\end{figure*}
Trajectories of the states of three agents under consensus protocol \eqref{eq:consensus_protocol} are depicted in Figure \ref{fig:sim_ex}.  Observe that $\m{x}_1^* = \m{x}_2^* = \m{x}_4^* \neq \m{x}_3^*$, as expected from the above discussion.
\end{Example}

Based on Proposition 1, we can further develop the following Corollaries.


\begin{Corollary} \label{cor:tree-expansion} Suppose a vertex $i$ connects to a positive tree $\mc{T} \subset \mc{G}$ via the edge set $\mc{S}_i = \{(i,j)$, $j \in \mc{N}_i \cap \mc{V}(\mc{T})\}$. If $\sum_{(i,j) \in \mc{S}} \m{A}_{ij}$ is positive definite, then under the consensus protocol \eqref{eq:consensus_protocol}, the equilibrium state of agent $i$ is the same with the equilibrium state of all agents in $\mc{T}$. 
\end{Corollary}

\begin{IEEEproof}
\begin{figure}[b]
  \centering
  \includegraphics[width=3.5cm]{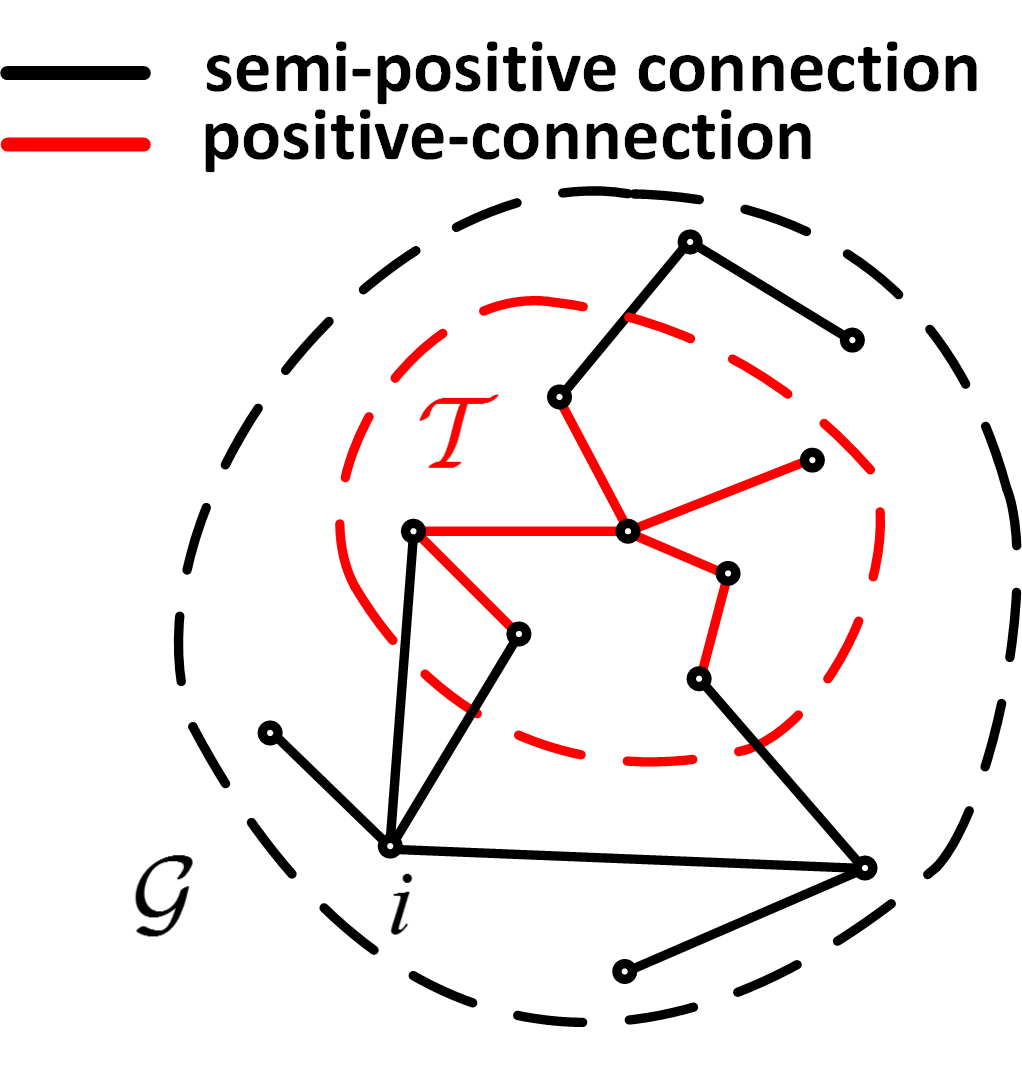}
  \caption{The vertex $i$ connects to the positive tree $\mc{T}$ through two semi-positive connections.}
  \label{fig:connect_tree}
\end{figure}
From Lemma \ref{lem:pos_tree}, we know that at equilibrium, all the states of all agents in the positive tree $\mc{T}$ are the same. Let $\m{x}^*_{\mc{T}}$ denote this value. Also, let $\m{x}_i^*$ denote the equilibrium value of agent $i$. For each semi-positive connection between $(i,j) \in \mc{S}$ (see Fig. \ref{fig:connect_tree} for an illustration), from Lemma \ref{lem:laplacian}, we have
\begin{equation*}
\m{A}_{ij} (\m{x}_i^* - \m{x}_{\mc{T}}^*) = \m{0}, \quad \forall (i,j) \in \mc{S}.
\end{equation*}
By adding the above equations, we can have
\begin{equation} \label{eq:expan-tree-1}
\left(\sum_{j \in \mc{S}} \m{A}_{ij} \right) (\m{x}_i^* - \m{x}_{\mc{T}}^*) = \m{0}.
\end{equation}
Since $\sum_{(i,j) \in \mc{S}} \m{A}_{ij}$ is positive definite, i.e. $\text{dim}(\bigcap_{j\in \mc{S}} \mc{N}(\m{A}_{ij}))=0$, equation  \eqref{eq:expan-tree-1} is satisfied if and only if $\m{x}_i^* = \m{x}_{\mc{T}}^*$.
\end{IEEEproof}

\begin{figure}[b]
\begin{center}
\includegraphics[width=3.6cm]{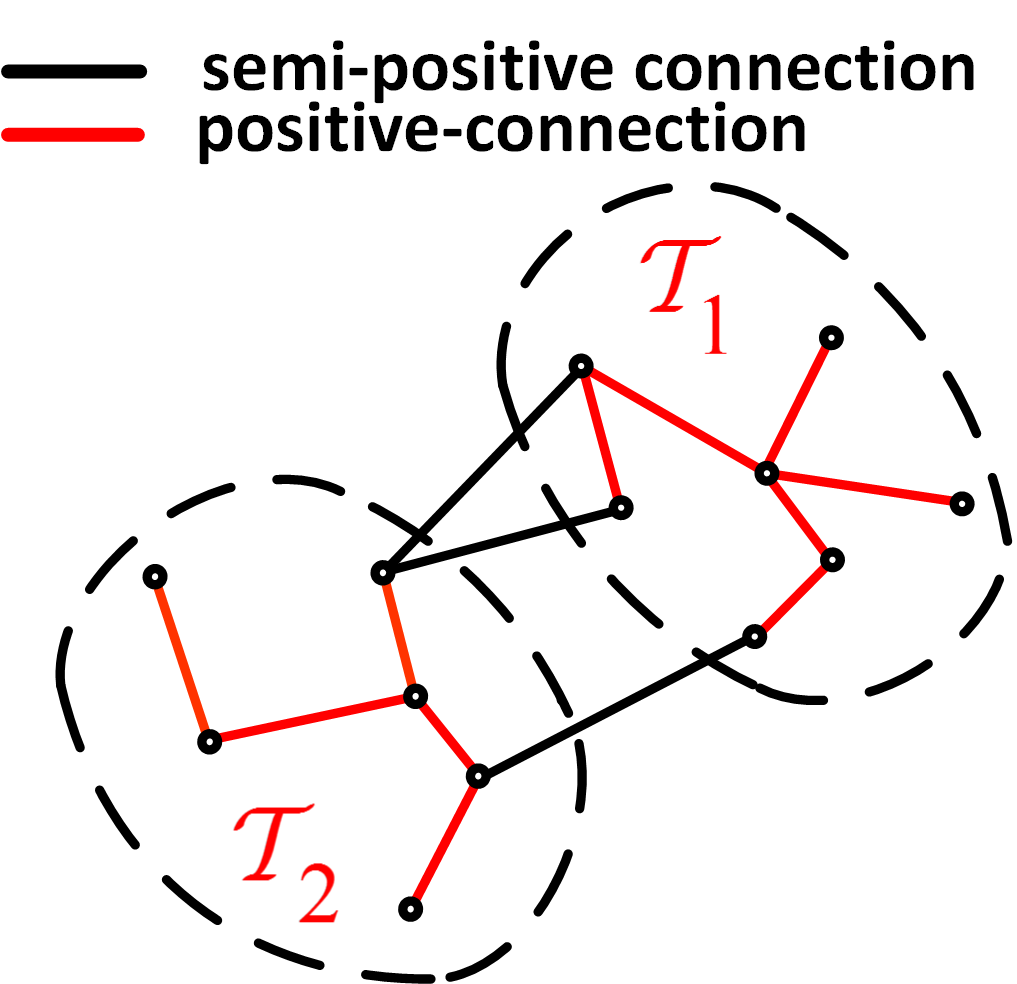}
\caption{Two positive trees $\mc{T}_1$ and $\mc{T}_2$ are connected through three semi-positive connections.}
\label{fig:two_trees}
\end{center}
\end{figure}
\begin{Corollary} \label{cor:EPT} Suppose two positive trees $\mc{T}_1, \mc{T}_2$ in  $\mc{G}$ are connected via the edge set $\mc{S}_i = \{ (i,j)|i \in \mc{T}_1, j \in \mc{T}_2 \}$. If $\sum_{(i,j) \in \mc{S}} \m{A}_{ij}$ is positive definite, then under the consensus protocol \eqref{eq:consensus_protocol}, the equilibrium states of all agents in $\mc{T}_1$ and $\mc{T}_2$ are the same.
\end{Corollary}

\begin{IEEEproof} Based on Lemma \ref{lem:pos_tree}, under \eqref{eq:consensus_protocol}, the equilibrium states of each agents belonging to the same positive tree are the same. Let $\m{x}_{\mc{T}_1}^*$ and $\m{x}_{\mc{T}_2}^*$ be the equilibrium states corresponding to each positive trees $\mc{T}_1$ and $\mc{T}_2$. For each semi-positive connection between $(i,j) \in \mc{S}$ (see Fig. \ref{fig:two_trees} for an illustration), from Lemma \ref{lem:laplacian}, we have
\begin{equation*}
\m{A}_{ij} (\m{x}_{\mc{T}_1}^* - \m{x}_{\mc{T}_2}^*) = \m{0}, \quad \forall (i,j) \in \mc{S}.
\end{equation*}
Adding the above equations, and following the same procedure as the proof of 
Corollary~\ref{cor:tree-expansion}, the proof can be completed.
\end{IEEEproof}

Figure~\ref{fig:two_trees} illustrates a scenario of Corollary \ref{cor:EPT}. Two positive trees $\mc{T}_1$ and $\mc{T}_2$ are connected through three semi-positive connections. If the summation of three matrices associated with these connections is positive definite, the equilibrium states of agents in both trees are the same.

\begin{Lemma}[Partitioning a graph into positive trees]\label{lem:partition}
Given a graph $\mc{G}$, consider a set of positive trees $\{\mc{T}_1, \ldots, \mc{T}_p \}$ ($1 \leq p \leq n$), where
\begin{itemize}
\item[(i)] $\mc{V}(\mc{T}_m) \bigcap \mc{V}(\mc{T}_l) = \emptyset$, $\bigcup_{m=1}^p \mc{V}(\mc{T}_m) = \mc{V}(\mc{G})$,
\item[(ii)] For each $\mc{T}_k$ ($1 \leq k \leq p$), $i, j \in \mc{V}(\mc{T}_k)$ if and only if there exists a positive path from $i$ to $j$.
\end{itemize}
Then, the partition of $\mc{G}$ defined by $\{\mc{V}(\mc{T}_1), \ldots, \mc{V}(\mc{T}_p)\}$ is unique.
\end{Lemma}

\begin{IEEEproof}
The following process gives a constructive way to find the partition:

We first select vertex 1. The tree $\mc{T}_1$ contains vertex 1 and all vertices that have a positive path to vertex 1 is unique. We then cross out all vertices in $\mc{T}_1$. The remaining vertices in $\mc{V}(\mc{G})\setminus \mc{V}(\mc{T}_1)$ do not have a positive path to any vertices in $\mc{V}(\mc{T}_1)$. Note $\mc{T}_1$ contains at least one vertex (vertex 1). Thus, $|\mc{V}(\mc{G})\setminus \mc{V}(\mc{T}_1)| < |\mc{V}(\mc{G})|$.

Next, we choose the vertex in $\mc{V}(\mc{G})\setminus \mc{V}(\mc{T}_1)$ with the smallest indexing. Similar to Step 1, we find the positive tree $\mc{T}_2$ associated with this vertex, and then cross out all vertices in $\mc{T}_2$ from the current vertex set. The remaining vertices are $\mc{V}(\mc{G})\setminus (\mc{V}(\mc{T}_1) \cup \mc{V}(\mc{T}_2))$, which has less vertices than  $\mc{V}(\mc{G})\setminus \mc{V}(\mc{T}_1)$.

We continue these processes, until there is no leftover vertex after crossing out all vertices from the last positive tree, say $\mc{T}_p$. At this point, we obtain a set of positive trees $\{\mc{T}_1, \ldots, \mc{T}_p \}$ ($1 \leq p \leq n$). Obviously, this set satisfies both conditions (i) and (ii).

Because in each step, the vertex and the corresponding positive tree are unique, the partition $\{\mc{V}(\mc{T}_1), \ldots, \mc{V}(\mc{T}_p)\}$ is unique.
\end{IEEEproof}

Let $\mc{C}(\mc{T}_m)$ be the cluster generated from the positive tree $\mc{T}_m$. If there exists a vertex $i \in \mc{C}(\mc{T}_m)$ satisfying the condition (ii) in Proposition \ref{prop:cluster} with a cluster $\mc{C}(\mc{T}_l)$, we can form a new cluster $\mc{C}(\mc{T}_m) \cup \mc{C}(\mc{T}_l)$ by merging $\mc{C}(\mc{T}_m)$ and $\mc{C}(\mc{T}_l)$ together. By this way, we can extend the positive trees in the graph. All vertices in the new cluster will reach a consensus under \eqref{eq:consensus_protocol}. To check whether two clusters $\mc{C}(\mc{T}_m)$ and $\mc{C}(\mc{T}_l)$ can be merged or not, it is sufficient to check condition (ii) in Proposition \ref{prop:cluster} for only one vertex $i \in \mc{C}(\mc{T}_m)$ with regard to $\mc{C}(\mc{T}_l)$. This property comes from the fact that \eqref{eq:same_cluster} is invariant for all vertices belonging to a same cluster.

Algorithm \ref{alg:cluster} proposes a solution for finding all clusters in the graph by iteratively checking condition (ii) in Proposition \ref{prop:cluster} and merging clusters together. The algorithm terminates after some finite steps and the output is a set of clusters $\mc{C}_{\mc{G}} = \{ \mc{C}_1, \ldots, \mc{C}_q \}$ ($1 \leq q \leq p$) satisfying
 \begin{itemize}
\item $\mc{C}_m \bigcap \mc{C}_l = \emptyset$, $\bigcup_{m=1}^q \mc{C}_m = \mc{V}(\mc{G})$,
\item For $1 \leq m \neq l \leq q$, $\mc{C}_m$ and $\mc{C}_l$ cannot be merged together, or i.e. $\nexists i \in \mc{C}_m$ satisfying the condition (ii) in Proposition \ref{prop:cluster} with $\mc{C}_l$.
\end{itemize}

We can now state the main result of this section.
\begin{Theorem} \label{thm:alg_graph_condition} Under the consensus protocol \eqref{eq:consensus}, the average consensus is achieved if and only if $\mc{G}$ is spanned by a cluster.
\end{Theorem}
\begin{IEEEproof}
If $\mc{G}$ is spanned by a cluster, it follows from Proposition \ref{prop:cluster} that the equilibrium state of all agents in the graph are the same, i.e., $\m{x}^*$ is the only equilibrium point of \eqref{eq:consensus}. Thus, the consensus is achieved globally exponentially  based on Lemma \ref{thm:stability}. On the other hand, if there is no cluster spanning $\mc{G}$, it follows that the agents belonging to two different clusters of $\mc{G}$ may not agree. Thus, a consensus cannot be globally achieved.
\end{IEEEproof}

\begin{algorithm}
\caption{Finding clusters of a matrix-weighted graph $\mc{G}$}
\label{alg:cluster}
\begin{algorithmic}[1]
\REQUIRE $\mc{G}(\mc{V},\mc{E},\mc{A})$
\STATE $i \leftarrow 0$;
\STATE Find the set of positive trees $\{ \mc{T}_1, \ldots, \mc{T}_p\}$ in $\mc{G}$;
\STATE $\mc{C}_\mc{G}(0) \leftarrow \{ \mc{C}_m = \{ \mc{V}(\mc{T}_m)\}, m = 1, \ldots, p \}$;
\REPEAT
	\STATE $\mc{C}_\mc{G}(i+1) \leftarrow \mc{C}_\mc{G}(i)$;
	\STATE check $\leftarrow$ false;
	\FORALL{$\mc{C}_m \in \mc{C}_\mc{G}(i)$}
		\FORALL{$\mc{C}_l \in \mc{C}_\mc{G}(i), l \neq m$}
    		\IF{$\exists i \in \mc{C}_l$ satisfies Theorem \ref{prop:cluster}(ii)}
    			\STATE $\mc{V}_{temp} \leftarrow \mc{V}(\mc{T}_m) \cup \mc{V}(\mc{T}_l)$;
    			\STATE $\mc{E}_{temp} \leftarrow \mc{E}(\mc{T}_m) \cup \mc{E}(\mc{T}_l) \cup \mc{S}$;
    			\STATE $\mc{C}_{temp} \leftarrow \mc{C}_m \cup \mc{C}_l$;
    			\STATE $\mc{C}_\mc{G}(i+1) \leftarrow (\mc{C}_\mc{G}(i+1) \setminus \{\mc{C}_m , \mc{C}_l\}) \cup \{\mc{C}_{temp}\}$;
    			\STATE check $\leftarrow$ true;
    			\STATE break;
   			\ENDIF
  		\ENDFOR
        \IF{check == true}
        \STATE break;
        \ENDIF
    \ENDFOR 
    \STATE $i$ $\leftarrow$ $i+1$; 
\UNTIL{$\mc{C}_\mc{G}(i) == \mc{C}_\mc{G}(i-1)$}
\end{algorithmic}
\end{algorithm}

\begin{Remark} Obviously, if the graph $\mc{G}$ has some disconnected components, under consensus protocol \eqref{eq:consensus_protocol}, the dynamics of each disconnected component do not influence the others. Suppose $\mc{G}$ is positive semiconnected and has $p$ clusters after Algorithm 1 terminates. Define $\mc{S}_{ij} := \{\mc{P}_k|~\text{the starting (end) vertex of $\mc{P}_k$ is in }~\mc{C}_i~(\text{resp.,} ~\mc{C}_j)\}$, the end states of each cluster satisfy:
\begin{align} 
\sum_{i=1}^q |C_i| \m{x}_{\mc{C}_i}^* &= n \bar{\m{x}}, \label{eq:constraints1}\\
\m{x}_{\mc{C}_i}^* - \m{x}_{\mc{C}_j}^* & \in \bigcap_{k=1}^{|\mc{S}_{ij}|} \mc{N}(\mc{P}_k).\label{eq:constraints2}
\end{align}
The solutions of equations \eqref{eq:constraints1}--\eqref{eq:constraints2} depend on matrix weights. Thus, we can design the matrix weights to obtain the desired number of clusters.
\end{Remark}

\begin{table*}[h!]
\centering
\caption{Comparison between the scalar-weighted consensus and the matrix-weighted consensus.} \label{table:1}
\begin{tabular}{|l|l|l|}
        \hline
        Property                                     & Scalar-weighted consensus                                & Matrix-weighted consensus                                                                                             \\ \hline
        Information flow                             & Fixed undirected connected                               & Fixed undirected semipositively connected                                                                             \\ 
        ~                                            & scalar-weighted graph: $\mc{G} = (\mc{V},\mc{E},\mc{A})$ & matrix-weighted graph: $\mc{G} = (\mc{V},\mc{E},\mc{A})$                                                              \\ \hline
        Edge's weights                               & $a_{ij}>0$: $(i,j)$ exists.                                & $\m{A}_{ij}>0$: $(i,j)$ is a positive definite edge.                                                                        \\ 
        ~                                            & $a_{ij}=0$: $(i,j)$ does not exist.                        & $\m{A}_{ij}\geq 0$: $(i,j)$ is a positive semidefinite edge.                                                                \\ 
        ~                                            & ~                                                        & $\m{A}_{ij} = 0$: $(i,j)$ does not exist.                                                                                   \\ \hline
        Graph Laplacian                            & $\m{L} = \m{H}^T diag(a_k) \m{H}$                        & $\m{L} = (\m{H}^T \otimes \m{I}_d) blkdiag(\m{A}_k) (\m{H} \otimes \m{I}_d)$                                                                                                
\\        ~                                            & $dim(\mc{N}(\m{L})) \geq 1$                              & $dim(\mc{N}(\m{L})) \geq d$                                                                                           \\ \hline
        Consensus protocol                           & $\dot{x}_i = \sum_{j \in \mc{N}_i} (x_j - x_i)$          & $\dot{\m{x}}_i = \sum_{j \in \mc{N}_i} \m{A}_{ij} (\m{x}_j - \m{x}_i)$                                                            \\ 
        ~                                            & $\dot{\m{x}} = -\m{L}\m{x}$                              & $\dot{\m{x}} = -\m{L}\m{x}$                                                                                           \\ \hline
        Consensus space                              & $span(\m{1}_n)$                                              & $span(\m{1}_n \otimes \m{I}_d)$                                                                                               \\ \hline
        Conditions for reaching                      & There is a tree spanning all vertices of  $\mc{G}$.      & There is a cluster spanning all vertices in $\mc{G}$.                                                               \\ 
        an average consensus                         & $\mc{N}(\m{L}) = span(\m{1}_n)$.                            & $\mc{N}(\m{L}) = span(\m{1}_n \otimes \m{I}_d)$.                                                                             \\ \hline
        Cluster consensus                            & Happens if and only if $\mc{G}$ is not connected.        & Happens if and only if $\mc{G}$ is not connected 
\\        & & or there does not exist a cluster spanning 
\\ & & all vertices of $\mc{G}$. \\
\hline
Average computational cost & $O(n \times d \times |\bar{\mc{N}_i}|)$ & $O(n \times d^2\times |\bar{\mc{N}_i}|)$ \\
\hline
    \end{tabular}
\end{table*}

The following result follows from Theorems \ref{thm:consensus_condition} and \ref{thm:alg_graph_condition}.
\begin{Theorem} Given a matrix-weighted graph $\mc{G}$, there is a cluster spanning all vertices of $\mc{G}$ if and only if $\mc{N}(\m{L}) = \mc{R}$.
\end{Theorem}

Before ending this section, we refer readers to Table~\ref{table:1}, which gives a comparison between the usual consensus algorithm and the matrix-weighted consensus algorithm proposed in this paper. Note that in Table \ref{table:1}, we find the computational cost as follow: Suppose an agent $i$ has an average number of neighbor $|\bar{\mc{N}_i}|$. Each matrix multiplication $\m{A}_{ij} \m{x}_i$ needs $d^2$ multiplications and $d(d-1)$ summations. Thus, in average, we needs $O(n \times d^2 \times |\bar{\mc{N}_i}|)$ calculations for each update.
\section{Applications}
\label{sec:5}
\subsection{Clustered consensus}
This subsection presents an example of designing a clustered consensus network based on the matrix-weighted consensus algorithm. 
Consider a system of nine autonomous agents in the plane. We would like to gather the agents into three clusters and then rendezvous to the system's average. Assume that the agents can sense the relative position with regard to its neighbor and there is a common reference frame for all agents in the system. The state of each agent is represented by a vector $\m{p}_i = [x_i,y_i]^T \in \mb{R}^2$, and the consensus protocol is explicitly written as 
$$\dot{\m{p}}_i = \sum_{j\in \mc{N}_i} \m{A}_{ij} (\m{p}_j - \m{p}_i),\quad \forall i = 1, \ldots, 9.$$ 

In case 1, we would like to gather the agents into three clusters. For this, the matrix weights are chosen as follows: 
$\m{A}_{12}=\begin{bmatrix} 2&0\\0&1\end{bmatrix}, $
$\m{A}_{13} = \begin{bmatrix} 2&3\\3&5\end{bmatrix}, $
$\m{A}_{47} = \begin{bmatrix} 0&0\\0&1\end{bmatrix}, $
$\m{A}_{14} = \begin{bmatrix} 0.75& {- 0.433}\\ {- 0.433} & 0.25\end{bmatrix}, $
$\m{A}_{17} = \begin{bmatrix} 0.75& { 0.433}\\ { 0.433} & 0.25\end{bmatrix}, $
$\m{A}_{45} = \begin{bmatrix} 1& { 0.5}\\ { 0.5} & 1\end{bmatrix}, $
$\m{A}_{46} = \begin{bmatrix} {0.9518}&{ - 0.2142}\\{ - 0.2142}&{0.0482}\end{bmatrix}, $
$\m{A}_{56} = \begin{bmatrix} {1}&{0}\\{0}&{0}\end{bmatrix}, $
$\m{A}_{78} = \begin{bmatrix} {3}&{2}\\{2}&{3}\end{bmatrix}, $ and
$\m{A}_{89} = \begin{bmatrix} {2}&{0}\\{0}&{2}\end{bmatrix}.$
Observe that $(1, 4)$, $(1, 7)$, $(4,7)$, and $(5,6)$ are semi-positive connections while other connections are positive. Thus, $\mc{G}$ has three clusters: $\mc{C}_1=\{1, 2, 3\}$, $\mc{C}_2=\{4, 5, 6\}$, $\mc{C}_3=\{7, 8, 9\}$. The equilibrium positions of three clusters satisfy: $\sum_{i=1}^3|\mc{C}_i|\m{p}_{\mc{C}_i}^*=\sum_{i=1}^9\m{p}_i(0)$, $\m{p}_{\mc{C}_1}^*- \m{p}_{\mc{C}_2}^* \in \mc{N}(\m{A}_{14})$, $\m{p}_{\mc{C}_1}^*- \m{p}_{\mc{C}_3}^* \in \mc{N}(\m{A}_{17})$, and $\m{p}_{\mc{C}_1}^*- \m{p}_{\mc{C}_3}^* \in \mc{N}(\m{A}_{47})$.

\begin{figure}[t!]
\begin{center}
\includegraphics[width=.5\linewidth]{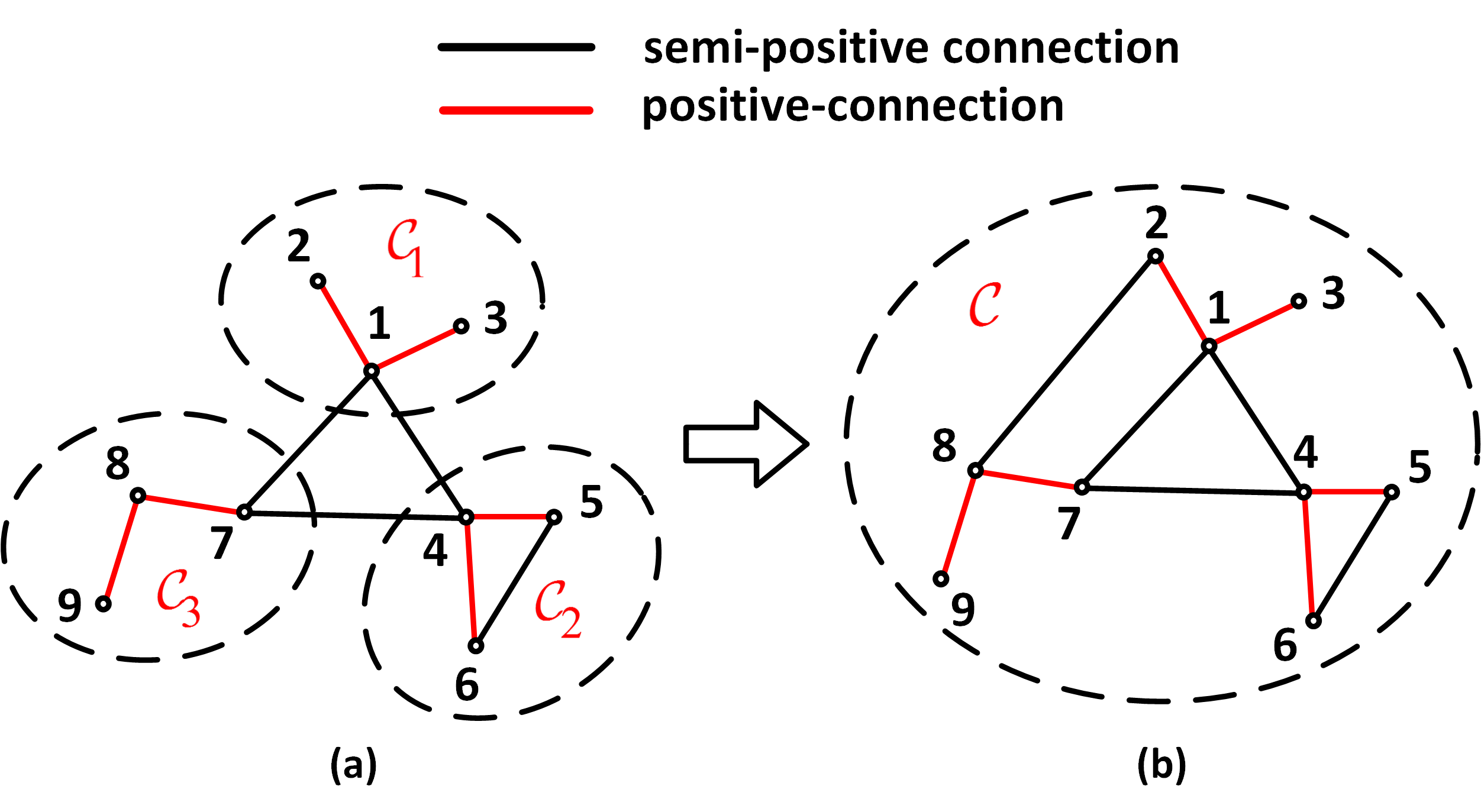}
\caption{Two graphs used in simulations: graph (a) has three clusters, graph (b) has a spanning cluster.}
\label{fig:cluster}
\end{center}
\end{figure}

In case 2, we want all agents to rendezvous at a point. To this end, we change the graph a little by adding a semipositive connection between vertices 2 and 8. The corresponding matrix weight is given by  $\m{A}_{28} = \begin{bmatrix} 0&0\\0&1\end{bmatrix}$. This additional connection makes two clusters $\mc{C}_1$ and $\mc{C}_3$ satisfy Corollary \ref{cor:EPT} (i.e., $(\m{A}_{17} + \m{A}_{28})$ is positive definite); thus they can be merged into a cluster, called $\mc{C}_{13}$. It is also due to Corollary \ref{cor:EPT} that the clusters $\mc{C}_{13}$ and $\mc{C}_{2}$ can be merged together (i.e. $(\m{A}_{14} + \m{A}_{17})$ is positive definite). It follows that the graph has a spanning cluster $\mc{C}$ in this case. Therefore, the agent will reach to an average consensus.

\begin{figure*}[t!]
    \begin{subfigure}[b]{0.33\textwidth}
    \centering
    \includegraphics[height=5cm]{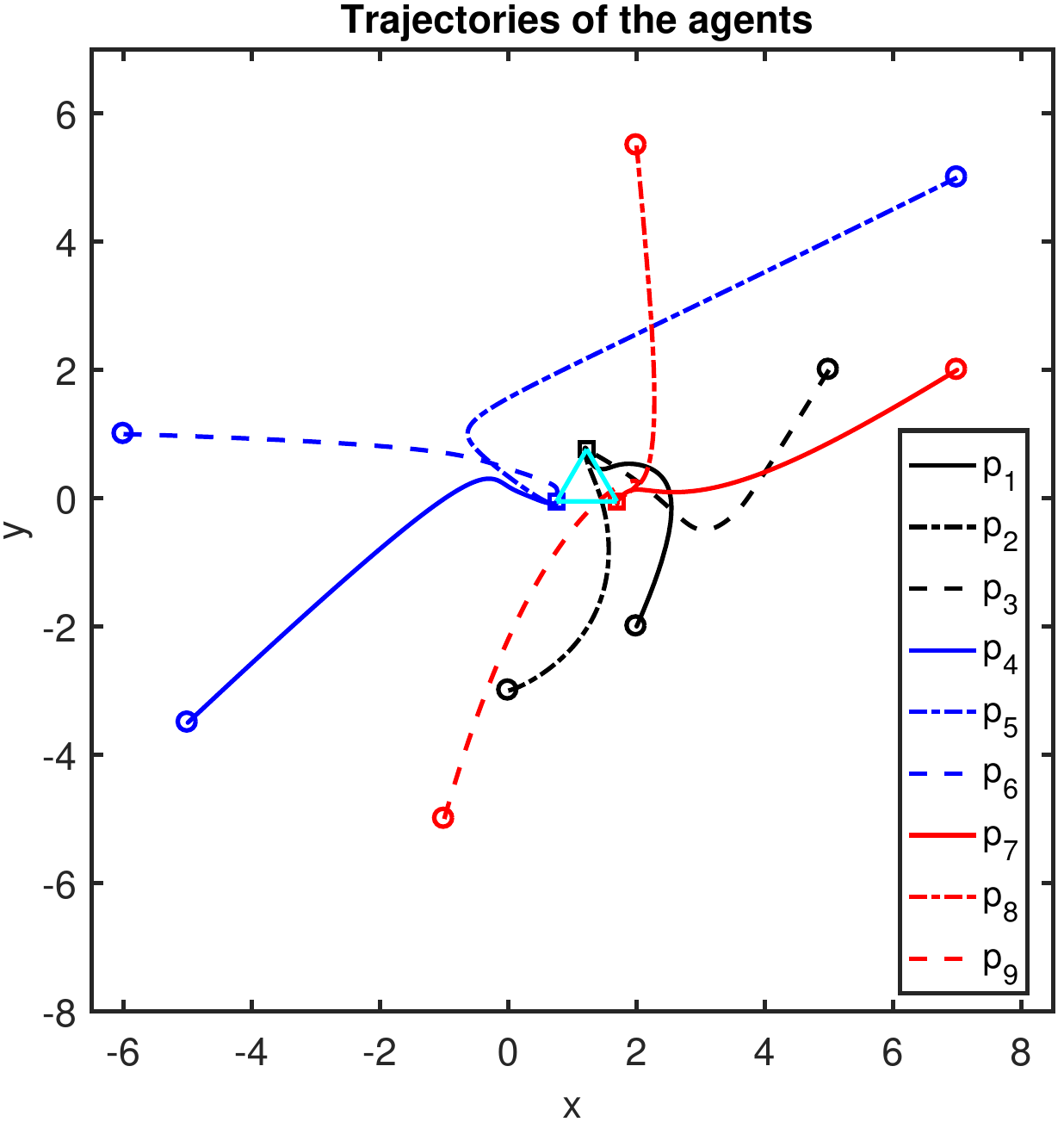}
    \caption{Trajectories of 9 agents.}
    \label{fig:sim1a}
    \end{subfigure}
    \begin{subfigure}[b]{0.33\textwidth}
    \centering
    \includegraphics[height=5cm]{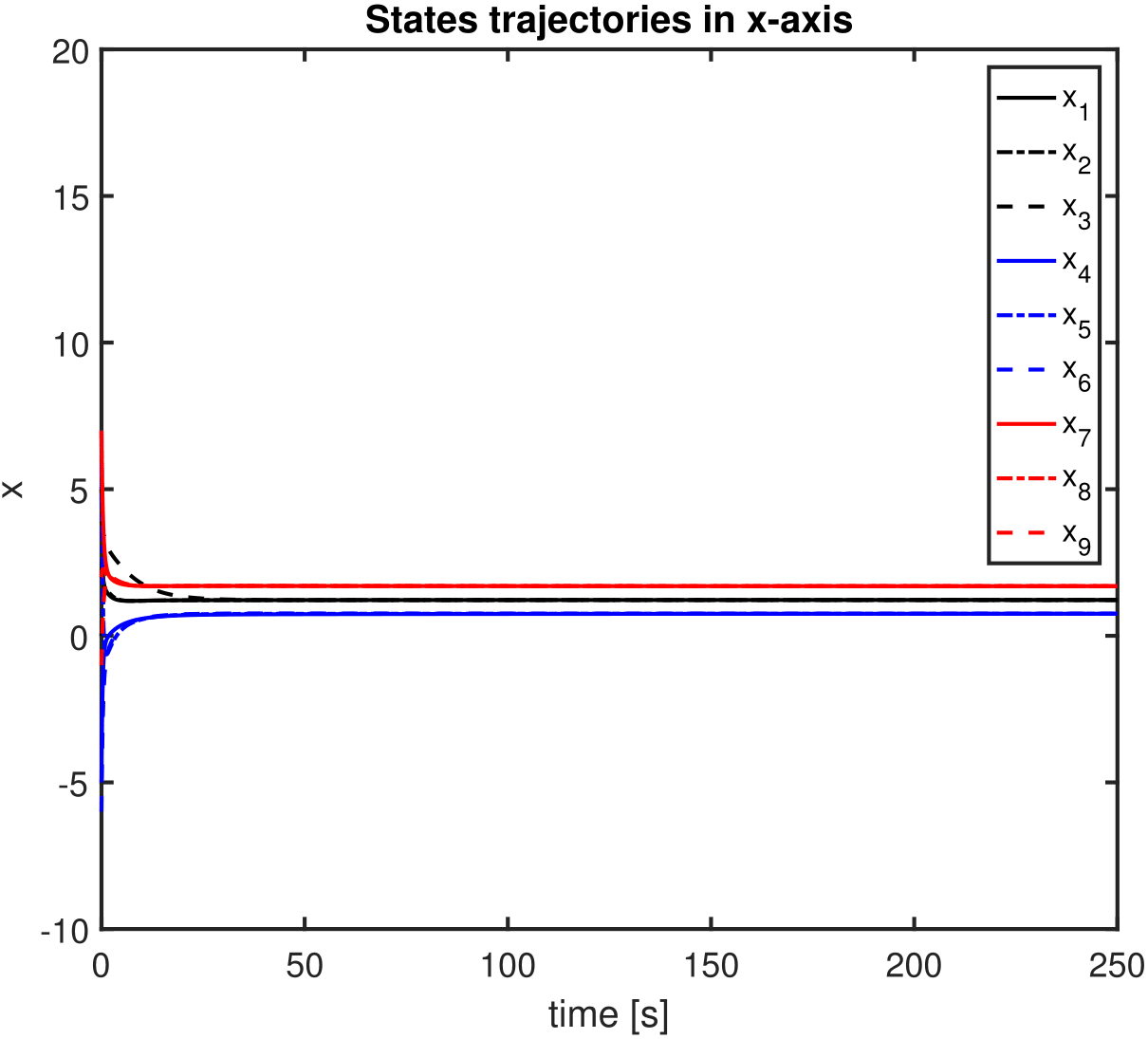}
    \caption{The $x$-axis dynamics.}
    \label{fig:sim1c}
    \end{subfigure}    
    \begin{subfigure}[b]{0.33\textwidth}
    \centering
    \includegraphics[height=5cm]{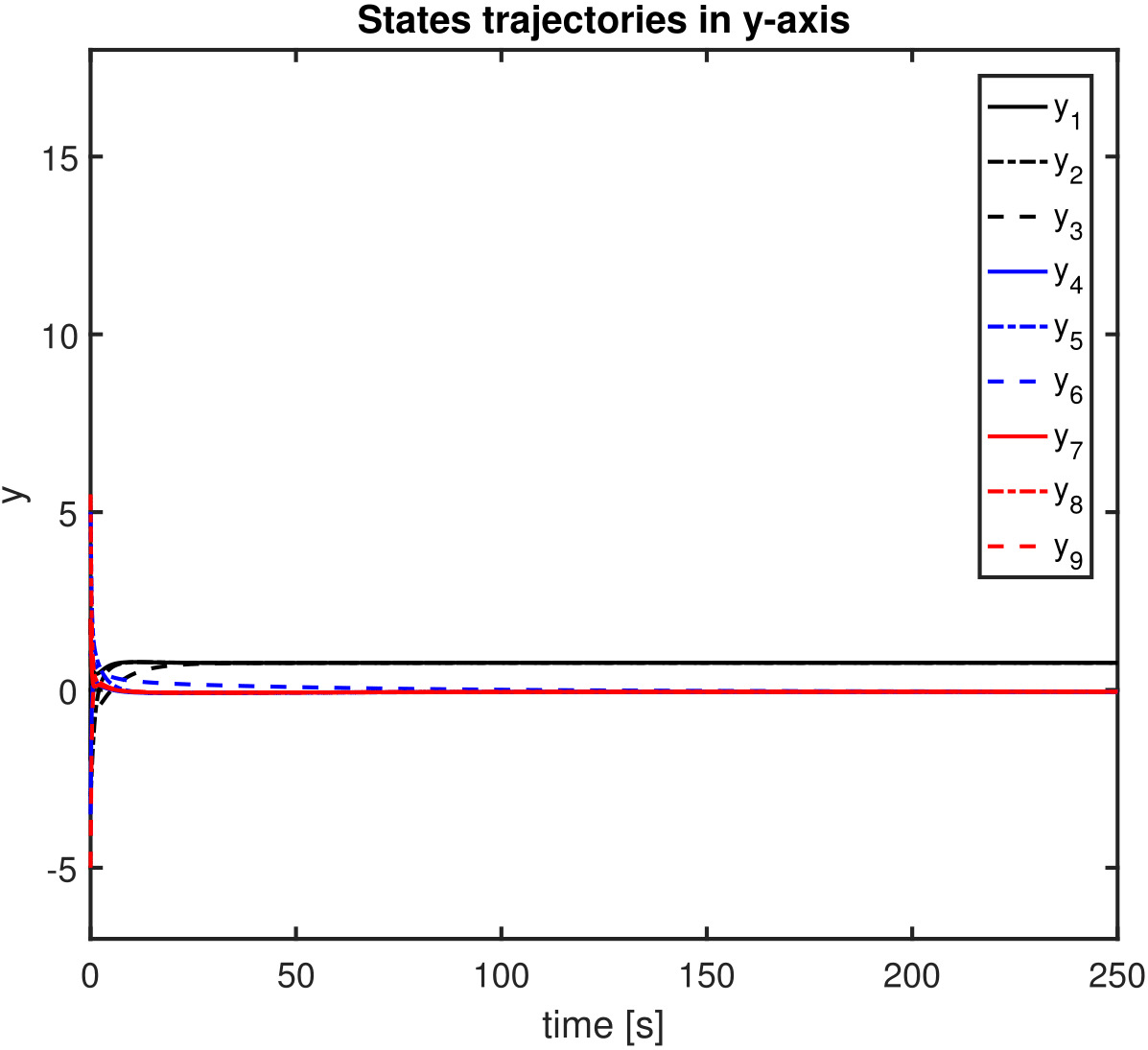}
    \caption{The $y$-axis dynamics.}
    \label{fig:sim1d}
    \end{subfigure}
    \caption{Case 1: The 9-agent system in Fig. \ref{fig:cluster} under the consensus protocol \eqref{eq:consensus_protocol}.}
    \label{fig:simulation}
\end{figure*}

We simulate the nine-agent system under the consensus protocol \eqref{eq:consensus_protocol}. Simulation results in case 1 are shown in Figure \ref{fig:simulation}. The state trajectories in $x$- and $y$-axes are depicted in Figures \ref{fig:sim1c} and \ref{fig:sim1d}, respectively. The corresponding positions of nine agents are shown in Figures \ref{fig:sim1a}. It can be seen that all agents belong to a cluster converge to a same point in $\mb{R}^2$.

\begin{figure*}[t]
    \begin{subfigure}[b]{0.33\textwidth}
    \centering
    \includegraphics[height=4.5cm]{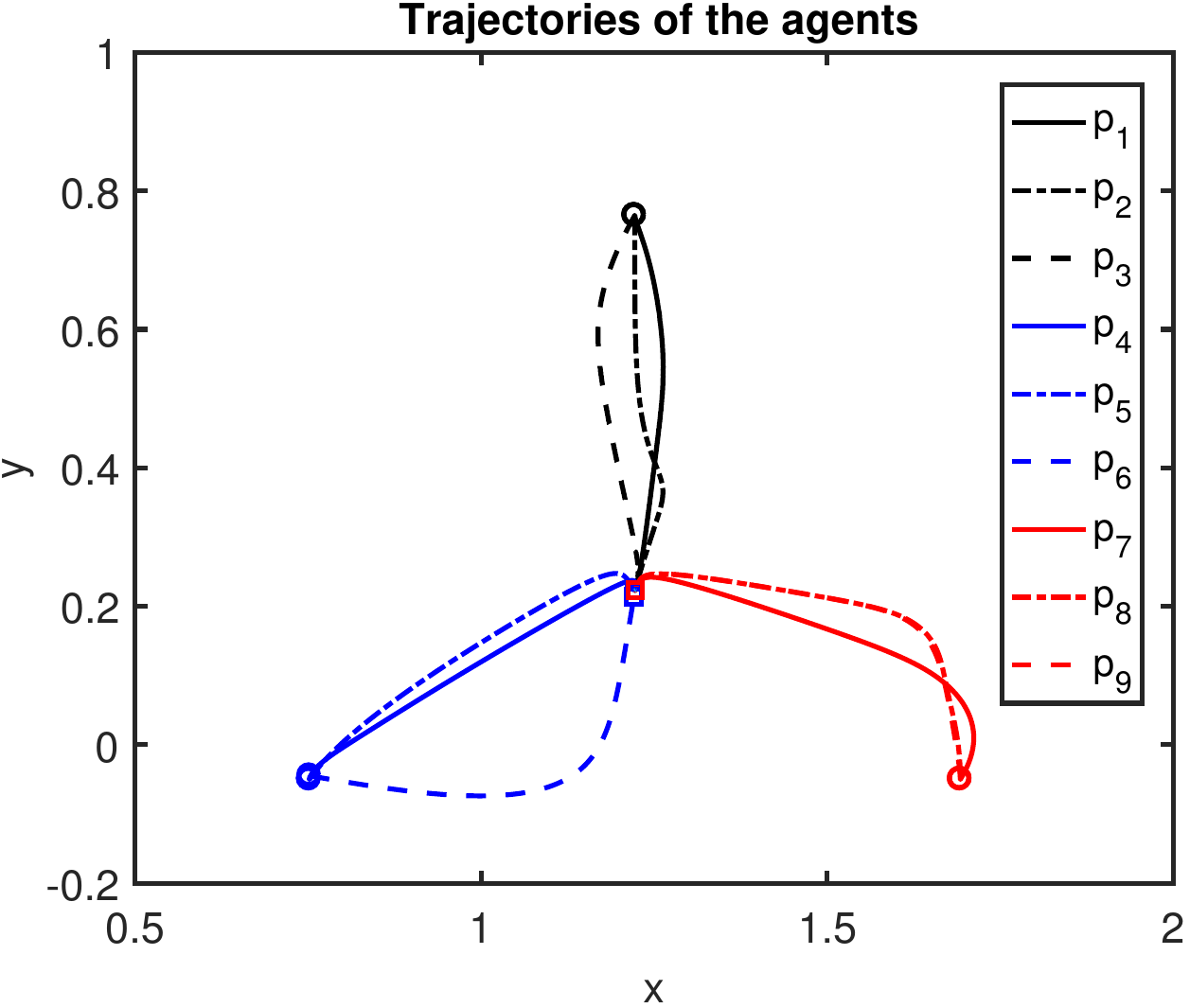}
    \caption{Trajectories of 9 agents.}
    \label{fig:sim1a1}
    \end{subfigure}
    \begin{subfigure}[b]{0.33\textwidth}
    \centering
    \includegraphics[height=4.5cm]{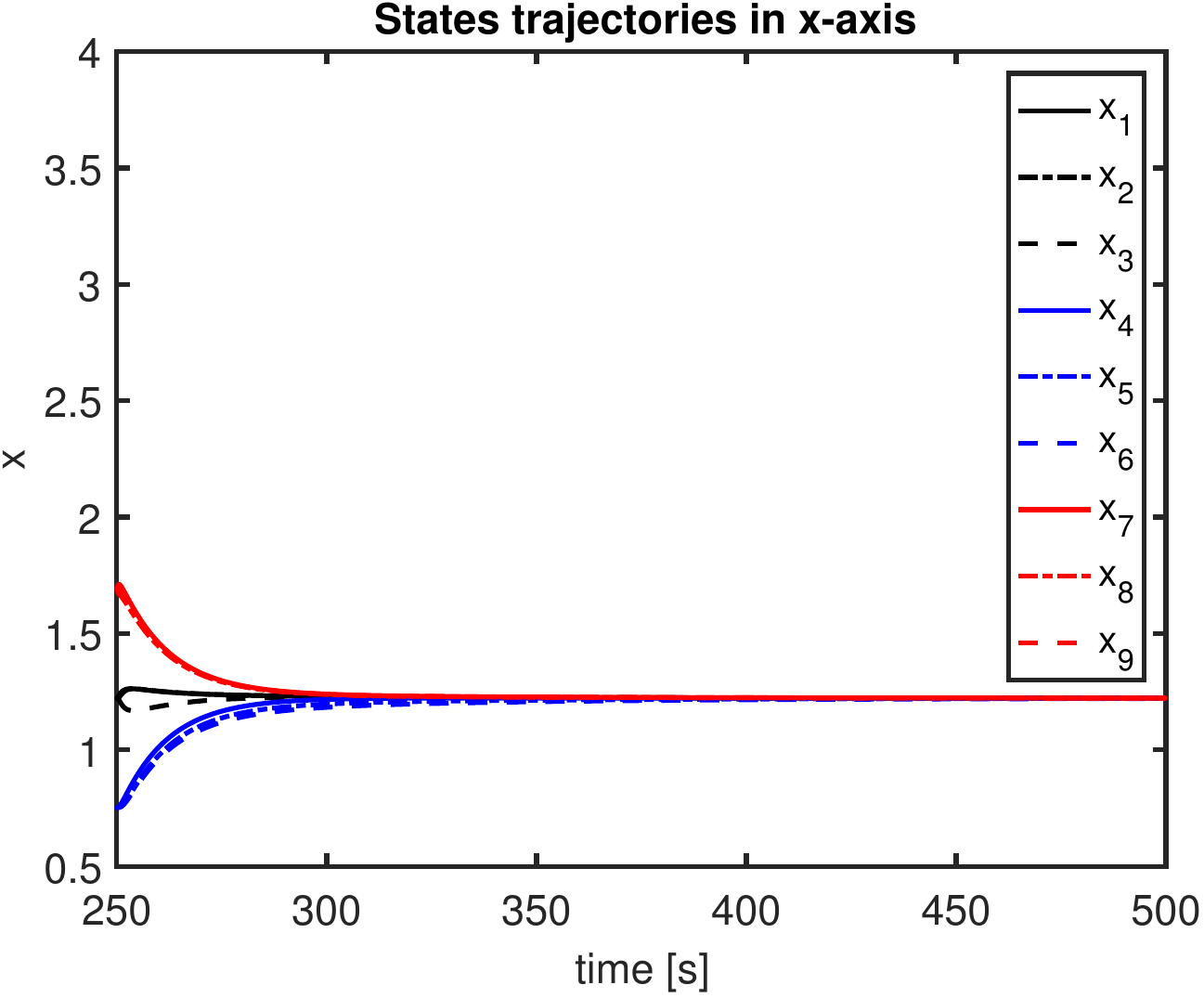}
    \caption{The $x$-axis dynamics.}
    \label{fig:sim1c1}
    \end{subfigure}    
    \begin{subfigure}[b]{0.33\textwidth}
    \centering
    \includegraphics[height=4.5cm]{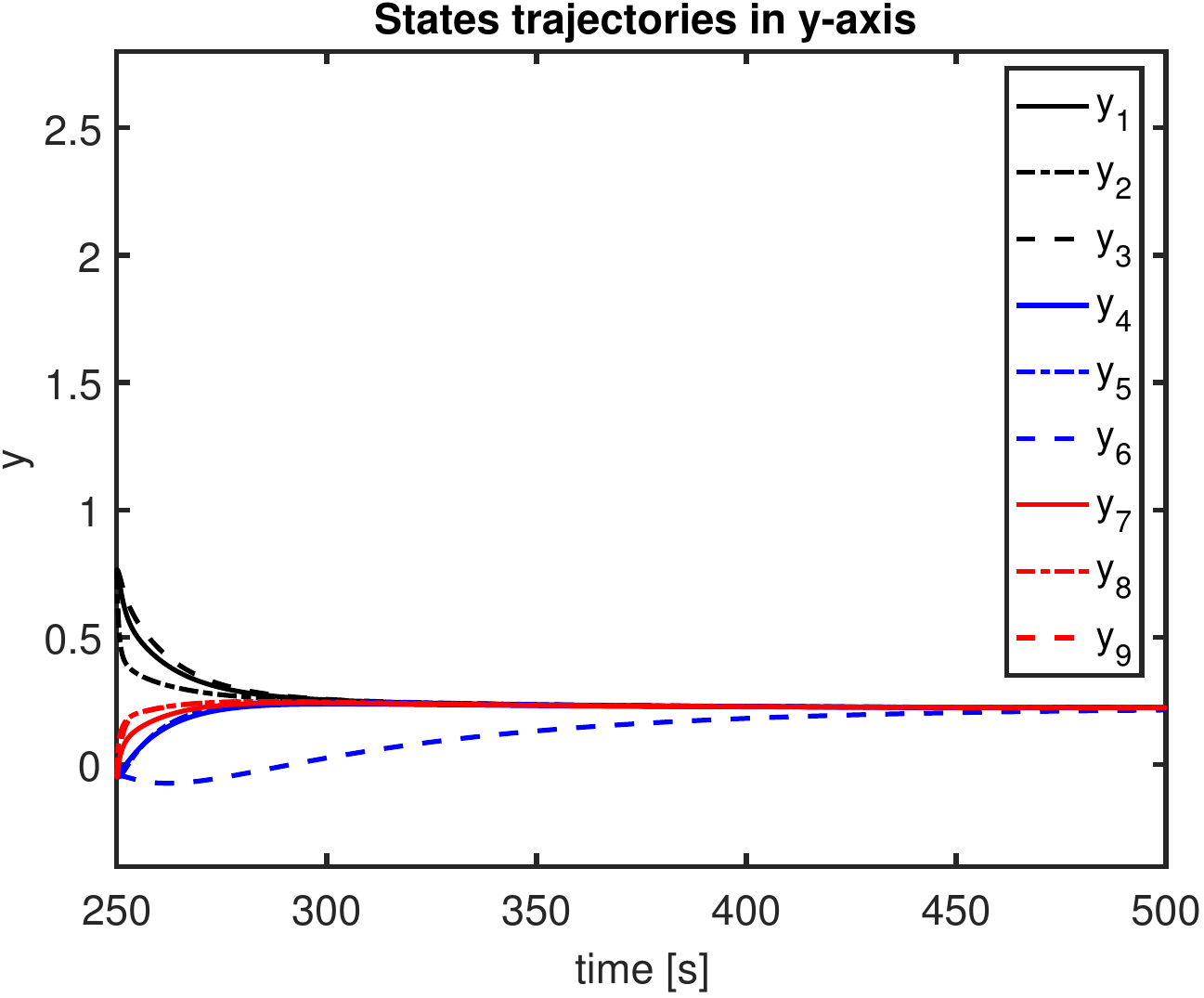}
    \caption{The $y$-axis dynamics.}
    \label{fig:sim1d1}
    \end{subfigure}
    \caption{Case 2: All agents rendezvous to a point after the interaction graph switched to the graph in Fig. \ref{fig:cluster} (ii).}
    \label{fig:simulation1}
\end{figure*}

Simulation results in case 2 are shown in Figure \ref{fig:simulation1}. Figure \ref{fig:sim1a1} depicts the nine agent trajectories after the interaction graph is switched to Fig. \ref{fig:cluster} (b). All agents asymptotically reach to a point in the plane.


\subsection{The bearing-constrained formation control problem}
The bearing-based formation control problem in \cite{Zhao2015CDC} can be considered as a special application of the matrix-weighted consensus problem. Here, the proposed formation control law for each agent is given by
\begin{equation} \label{eq:bearing_based_formation_control}
\dot{\m{p}}_i = - \sum_{j \in \mc{N}_i} \m{P}_{\m{g}_{ij}^*}(\m{p}_i - \m{p}_{j}),
\end{equation}
where $\m{p}_i \in \mb{R}^d$ is the position of agent $i$. The control law \eqref{eq:bearing_based_formation_control} is exactly the consensus protocol \eqref{eq:consensus_protocol} where the matrix weights are chosen to be the projection matrices $\m{P}_{\m{g}_{ij}^*}$. Here, $\m{g}_{ij}^*$ is a unit bearing vector which has been chosen to impose a constraint for the formation. Also, the projection matrix is defined as  $\m{P}_{\m{g}_{ij}^*}:= \m{I}_{d \times d} - \m{g}_{ij}^* \m{g}_{ij}^{*T}$, and thus it is symmetric, positive semidefinite. Moreover, the nullspace of $\m{P}_{\m{g}_{ij}^*}$ is spanned by the bearing vector $\m{g}_{ij}^*$, i.e., $\mc{N}(\m{P}_{\m{g}_{ij}^*}) = \text{span}\{\m{g}_{ij}^*\}$.

By specifying a set of desired bearing vectors $\{\m{g}_{ij}^*\}_{(i,j) \in \mc{E}}$ to a desired formation $\m{p}^*$, we can design the nullspace of the \emph{bearing Laplacian matrix} $\m{L}_B(\m{p}^*) \in \mb{R}^{dn\times dn}$.
Note that the $ij$th block sub-matrix of $\m{L}_B$ is given by
\[\left\{ {\begin{array}{*{20}{l}}
  {{{\left[ {{\m{L}_B}} \right]}_{ij}} = \m{0},}&{i \ne j,\left( {i,j} \right) \notin \mc{E},} \\ 
  {{{\left[ {{\m{L}_B}} \right]}_{ij}}=  - {\m{P}_{\m{g}_{ij}^*}},}&{i \ne j,\left( {i,j} \right) \in \mc{E},} \\ 
  {\left[ {{\m{L}_B}} \right]_{ii} = \sum\nolimits_{j \in {\mc{N}_i}} {{\m{P}_{\m{g}_{ij}^*}},} }&{i \in \mc{V}.} 
\end{array}} \right.\]

Based on Theorem \ref{thm:stability}, the formation converges to the nullspace of $\m{L}_B$ which is $\text{ span }\{ \mc{R}, \{\m{v} = [\m{v}_1^T, \ldots, \m{v}_n^T]^T \in \mb{R}^{dn}| (\m{v}_j - \m{v}_i) = \alpha \m{g}_{ij}^*), \alpha \in \mb{R}, \forall (i,j) \in \mc{E} \}$.
More detailed analysis and discussions on the bearing-based formation control can be found in \cite{Zhao2015CDC,Zhao2015CNS,Zhao2016aut}.

\section{Conclusion}
\label{sec:6}
In this paper, the matrix-weighted consensus algorithm was proposed. It was shown that the matrix-weighted consensus algorithm exhibits both common and unique characteristics compared with the usual consensus algorithm in \cite{Olfati2004}. More specifically, under the matrix-weighted consensus algorithm, connectedness of the undirected graph is not enough to guarantee the system to globally achieve a consensus. In fact, due to the existence of semipositive connections, the clustered consensus phenomena can easily happen in the network. It was proved that a global consensus can be achieved if and only if the network is spanned by only one cluster. Further, an algorithm for finding all clusters in the network was also provided. We illustrated two possible applications of the matrix-weighted consensus protocol in clustered consensus and in bearing-based formation control problem. 

There are still several open problems on matrix-weighted consensus. An immediate problem is finding a necessary and sufficient condition for achieving a consensus in directed networks. For this, matrix-weighted consensus with leader-following graphs have been recently studied in \cite{Trinh2017ASCC}. Examining the consensus algorithm with communication delays, or time-varying matrix weights, or switching topologies along the works \cite{cao2012distributed,kan2016leader}, can be several further research directions. Finally, the matrix-weighted consensus protocol can have applications in cyber-physical systems and in modeling of social networks.

%


\bibliographystyle{IEEEtran}
\bibliography{minh2017}

\end{document}